\documentclass[11pt,reqno]{amsart}      
\usepackage[english]{babel}
\usepackage[utf8x]{inputenc}
\usepackage[T1]{fontenc}
\usepackage[letterpaper, margin=1.2in]{geometry}
\usepackage{mathtools}
\usepackage[dvipsnames]{xcolor}
\usepackage{hyperref}
\usepackage{dsfont}
\usepackage{enumerate}
\usepackage{dsfont}
\usepackage[pdftex]{graphicx}
\usepackage{caption}
\usepackage{subcaption}
\usepackage{tikz}
\usetikzlibrary{positioning,decorations.pathreplacing,shapes.misc}
\tikzset{cross/.style={cross out, draw=black, fill=none, minimum size=2*(#1-\pgflinewidth), inner sep=0pt, outer sep=0pt}, cross/.default={2pt}}

\allowdisplaybreaks

\usepackage{amsmath, amsthm, amssymb}
\usepackage{latexsym,graphicx,multicol,mathrsfs,xypic}

\newtheorem{theorem}{Theorem}[section]
\newtheorem{definition}[theorem]{Definition}
\newtheorem{proposition}[theorem]{Proposition}
\newtheorem{lemma}[theorem]{Lemma}

\newtheorem{remark}[theorem]{Remark}

\newtheorem*{proposition*}{Proposition}

\DeclareMathOperator{\var}{var}
\DeclareMathOperator{\cov}{cov}

\DeclareMathOperator{\Hess}{Hess}
\DeclareMathOperator{\Id}{Id}
\DeclareMathOperator{\supp}{supp}

\usepackage{autonum}
\begin{document}

\title[Uniform LSI for the canonical ensemble]{Uniform LSI for the canonical ensemble on the 1d-lattice with strong, finite-range interaction.}

\author{Younghak Kwon}
\address{Department of Mathematics, University of California, Los Angeles}
\email{yhkwon@math.ucla.edu}

\author{Georg Menz}
\address{Department of Mathematics, University of California, Los Angeles}
\email{gmenz@math.ucla.edu}


\date{\today}

\begin{abstract}

We consider a one-dimensional lattice system of unbounded, real-valued spins with arbitrary strong, quadratic, finite-range interaction. We show that the canonical ensemble (ce) satisfies a uniform logarithmic Sobolev inequality (LSI). The LSI constant is uniform in the boundary data, the external field and scales optimally in the system size. This extends a classical result of H.T.~Yau from discrete to unbounded, real-valued spins. It also extends prior results of Landim, Panizo \& Yau or Menz for unbounded, real-valued spins from absent- or weak- to strong-interaction. We deduce the LSI by combining two competing methods, the two-scale approach and the Zegarlinski method. Main ingredients are the strict convexity of the coarse-grained Hamiltonian, the equivalence of ensembles and the decay of correlations in the ce.   
\end{abstract}

\maketitle

\section*{Introduction} \label{intro}
In this article, we study a one-dimensional lattice system of unbounded real-valued spins. The system consists of a finite number of sites~$i \in \Lambda \subset \mathbb{Z}$ on the lattice~$\mathbb{Z}$. For convenience, we assume that the set~$\Lambda$ is given by~$\left\{1, \ldots, N \right\}$. At each site~$i \in \Lambda$ there is a spin~$x_i$. In the Ising model the spins can take on the value~$0$ or~$1$. Here, spins~$x_i \in \mathbb{R}$ are real-valued and unbounded. A configuration of the lattice system is given by a vector~$x \in \mathbb{R}^{N}$. The energy of a configuration~$x$ is given by the Hamiltonian~$H: \mathbb{R}^N \to \mathbb{R}$ of the system. We consider arbitrary strong, pairwise, quadratic, finite-range interaction. For the detailed definition of the Hamiltonian~$H$ we refer to Section~\ref{s_setting_and_main_results}.\\

We consider two different ensembles: The first ensemble is the grand canonical ensemble (gce) which is given by the finite-volume Gibbs measure
\begin{align}
  \mu^\sigma(dx) = \frac{1}{Z} \exp\left( \sigma \sum_{i=1}^N x_i- H(x) \right) dx.
\end{align}
Here,~$Z$ is a generic normalization constant making the measure~$\mu^\sigma$ a probability measure. The constant~$\sigma \in \mathbb{R}$ is interpreted as an external field. The second ensemble is the canonical ensemble (ce). It emerges from the gce by conditioning on the mean spin
\begin{align} \label{e_mean_spin}
  m = \frac{1}{N} \sum_{i =1}^N x_i.
\end{align}
The ce is given by the probability measure
\begin{align}%\label{d_c_ensemble}
  \mu_m(dx)&= \mu^{\sigma} \left( dx \ | \    \frac{1}{N} \sum_{i=1}^N x_i= m  \right)\\
  &= \frac{1}{Z} \  \mathds{1}_{\left\{ \frac{1}{N}\sum_{i=1}^N x_i= m \right\}  } (x)\ \exp\left(\sigma \sum_{i=1}^{N} x_i - H(x)\right) \mathcal{L}^{N-1}(dx) \label{e_ce_cancellation}\\
  &= \frac{1}{Z} \  \mathds{1}_{\left\{ \frac{1}{N}\sum_{i=1}^N x_i= m \right\}  } (x)\ \exp(- H(x)) \mathcal{L}^{N-1}(dx), \label{e_ce_cancellation_result}
\end{align} 
where~$\mathcal{L}^{N-1}$ denotes the~$N-1$-dimensional Hausdorff measure. The last identity follows from the observation that $\sum_{i=1}^{N} x_i = Nm$ and therefore the factor~$\sigma \sum_{i=1}^N x_i = \sigma N m$ can be cancelled out by the same factor appearing in the normalization constant~$Z$. \\

The logarithmic Sobolev inequality (LSI) – introduced by Gross~\cite{Gro75}– is a powerful tool for studying spin systems. For example, the LSI implies Gaussian concentration via Herbst's argument, is equivalent to hypercontractivity and characterizes the exponential rate of convergence to equilibrium of the naturally associated diffusion process. By the equivalence of dynamic and static phase transitions, a uniform LSI also indicates the absence of a phase transition (see e.g.~\cite{Yos03,HeMe16}). For an introduction to the LSI we refer the reader to~\cite{Led01,Led01a,Roy99,BaGeLe14}. \\

On the one-dimensional lattice, a uniform LSI holds for the gce~$
\mu^\sigma$ even for infinite-range interactions, given the interaction decays fast enough (see~\cite{MeNi14}). Deducing a uniform LSI becomes a lot harder when considering the ce~$\mu_m$ instead of the gce~$\mu^\sigma$. Even if there is no interaction term in the Hamiltonian~$H$, the ce~$\mu_m$ is not a product measure. There are long-range interactions due to the conditioning onto the mean spin~$m= \frac{1}{N} \sum_{i=1}^N x_i$. For~$\left\{0,1 \right\}$-valued spins on arbitrary lattices, a classical result of Yau~\cite{Yau96} states that the ce~$\mu_m$ satisfies a uniform LSI as soon as the correlations of the gce~$\mu^{\sigma}$ decay exponentially, which is the case on the one-dimensional lattice. The original proof by Yau~\cite{Yau96} is based on the Lu-Yau Martingale Method~\cite{LuYa93}. Later, Cancrini, Martinelli \& Roberto~\cite{CMR02} gave an alternative, self-contained proof of the same statement.\\ 

Extending Yau's result~\cite{Yau96} to unbounded, real-valued spins is still an open problem. The main result of this article (cf.~Theorem~\ref{p_uniform_lsi} below) solves this problem. Considering unbounded, real-valued spins instead of $\left\{0,1 \right\}$-valued spins yields a technical challenge: Because spins are unbounded compactness is lost and many arguments do not carry over from the discrete case. Therefore, it was already quite a challenge to establish the LSI for the canonical ensemble in the case of a non-interacting Hamiltonian (cf.~\cite{Cha03,GrOtViWe09,LaPaYa02,MeOt13} for unbounded, real-valued spins). The main difficulty was to obtain the optimal scaling behavior of the LSI constant in the system size. 
More recently, the uniform LSI was deduced for arbitrary weak interaction in~\cite{Me11}. The method used in~\cite{Me11} is of perturbative nature and different ideas are needed when considering strong interaction. Therefore it is not surprising that deducing the LSI for the ce~$\mu_m$ under a mixing condition for the gce~$\mu^\sigma$ remained an open problem.\\

A major breakthrough was recently accomplished in~\cite{KwMe17} where a non-perturbative result on the equivalence of ensembles of the gce and the ce was deduced. In~\cite{KwMe17}, only attractive, nearest-neighbor interaction was considered. In Section~\ref{s_auxiliary_lemmas} we extend those results to arbitrary finite-range interaction. A consequence is that the coarse-grained Hamiltonian of a single block is uniformly strictly convex (see Corollary~2 in~\cite{KwMe17}). This rises hope that the uniform LSI for the ce could be deduced via the two-scale approach~\cite{GrOtViWe09}. It turns out that the situation is more complicated and the result of~\cite{KwMe17} is not sufficient to directly apply the two-scale approach. The strict convexity of the coarse-grained Hamiltonian was deduced in~\cite{KwMe17} for only one block. One would need strict convexity for all blocks simultaneously. This is a lot harder to show due to the strong interaction between blocks. However, we still manage to build up on the results of~\cite{KwMe17} and deduce the uniform LSI for ce under a mixing condition for the gce (see Theorem~\ref{p_uniform_lsi}). We overcome the obstacle of strong interactions between blocks by combining several ideas and methods from the two-scale approach (see~\cite{OttRez07,GrOtViWe09,Me13}), the Zegarlinski method (see~\cite{Zeg96}), decay of correlations (cf.~\cite{KwMe18a}) and a decomposition method for Hamiltonians introduced in~\cite{Me13}. For more details on the argument we refer to Section~\ref{s_setting_and_main_results} and Section~\ref{s_proof_uniform_lsi}.\\

Deducing a uniform LSI for the ce~$\mu_m$ has another special importance: It is one of the main ingredients when deducing the hydrodynamic limit of the Kawasaki dynamics via the two-scale approach~\cite{GrOtViWe09}. Because the uniform LSI controls the entropy production, it also plays an implicit role in other approaches to the hydrodynamic limit via the entropy method, the martingale method or the gradient flow method (see for example \cite{GuPaVa88,Yau91,KipLan99,DuoFat16}). The Kawasaki dynamics is a natural drift diffusion process on the lattice system that conserves the mean spin of the system. The ce~$\mu_m$ is the stationary and ergodic distribution of the Kawasaki dynamics. The hydrodynamic limit is a dynamic manifestation of the law of large numbers. It states that under the correct scaling the Kawasaki dynamics (which is a stochastic process) converges to the solution of a non-linear heat equation (which is deterministic). It is conjectured by H.T.~Yau that the hydrodynamic limit also holds for strong finite-range interactions on a one-dimensional lattice. So far, this conjecture also is wide open. Deducing a uniform LSI for the ce~$\mu_m$ on the one-dimensional lattice with arbitrary strong finite-range interaction is an important interim result to attack this problem.\\

Let us comment on open questions and problems:
\begin{itemize}

\item Instead for finite-range interaction, could one deduce similar results for infinite-range, algebraically decaying interactions? More precisely, is it possible to extend the results of~\cite{MeNi14} from the gce to the ce? If yes, is the same order of algebraic decay sufficient, i.e.~of the order~$2+ \varepsilon$, or does one need a higher order of decay? For solving this problem one would have to overcome several difficulties. For example, generalizing the equivalence of ensembles (see~\cite{KwMe17}) would need new work. Also, because we use ideas of the Zegarlinski method, the arguments of this article are restricted to the one-dimensional lattice with finite-range interaction. Applying our method to infinite-range interaction would yield a cyclic dependence of the different parameters. A possible alternative approach to this problem is to generalize the approach of~\cite{OttRez07,Me13,MeNi14} from the ce to the gce.\\[-1ex]

\item Can one show that, as it is the case for the gce, there is a phase transition for sufficiently slow decaying, infinite-range interaction (see for example~\cite{Dys69,FrSp82,CaFeMePr05,Imb82} for related results on the Ising model)?\\[-1ex]

\item Is it possible to consider more general Hamiltonians? For example, our argument is based on the fact that the single-site potentials are perturbed quadratic, especially when we use the results of~\cite{KwMe17}. One would like to have general super-quadratic potentials as was for example used in~\cite{MeOt13}. Also, it would be nice to consider general interactions than quadratic or pairwise interaction.\\[-1ex]

\item Is it possible to generalize the results to vector-valued spin systems?
\end{itemize}

\noindent There are many ways to proceed from this article: 

\begin{itemize}
\item Motivated by the results for discrete-spin systems (cf.~\cite{Yau96,CMR02}), one could try to deduce the hydrodynamic limit of the Kawasaki dynamics under a mixing condition for the gce on arbitrary lattices. The first strategy would be using the two-scale approach~\cite{GrOtViWe09}. However, this approach needs more non-trivial ingredients than the ones provided in this work and~\cite{KwMe17,KwMe18a}.\\[-1ex]

\item Is for the ce~$\mu_m$ a uniform LSI equivalent to decay of correlations, as it is the case for the gce~$\mu^\sigma$ (see for example~\cite{Yos01,Yos03,HeMe16})?\\[-1ex]

\end{itemize}

We conclude the introduction by giving an overview over the article. In Section~\ref{s_setting_and_main_results} we introduce the precise setting and present the main results. In Section~\ref{s_auxiliary_lemmas} we provide several auxiliary results. In Section~\ref{s_proof_uniform_lsi} we give the proof of the main result of this article, namely the uniform LSI for the ce (see Theorem~\ref{s_proof_uniform_lsi}). The proof of Theorem~\ref{s_proof_uniform_lsi} is based on three auxiliary statements, which are deduced in Section~\ref{s_proof_conditional_lsi}, Section~\ref{s_proof_marginal_lsi} and Section~\ref{s_proof_combined_lsi}.

\section*{Conventions and Notation}

\begin{itemize}
\item The symbol~$T_{(k)}$ denotes the term that is given by the line~$(k)$.
\item We denote with~$0<C<\infty$ a generic uniform constant. This means that the actual value of~$C$ might change from line to line or even within a line.
\item Uniform means that a statement holds uniformly in the system size~$N$, the mean spin~$m$ and the external field~$s$.
\item $a \lesssim b$ denotes that there is a uniform constant~$C$ such that~$a \leq C b$.
\item $a \sim b$ means that~$a \lesssim b$ and~$b \lesssim a$.
\item $\mathcal{L}^{k}$ denotes the $k$-dimensional Hausdorff measure. If there is no cause of confusion we write~$\mathcal{L}$.
\item $Z$ is a generic normalization constant. It denotes the partition function of a measure.  
\item For each~$N \in \mathbb{N}$,~$[N]$ denotes the set~$\left\{ 1, \ldots N \right\}$.
\item For a vector~$x \in \mathbb{R}^{N}$ and a set~$A \subset [N]$,~$x^A \in \mathbb{R}^{A}$ denotes the vector $ (x^A)_{i} = x_i$ for all~$i \in A$.
\item For a function~$f : \mathbb{R}^{N} \to \mathbb{C}$, we denote with~$\supp f = \{i_1, \cdots, i_k \}$ the minimal subset of~$[N]$ such that~$f(x) = f(x_{i_1}, \cdots, x_{i_k})$. 
\end{itemize}

\section{Setting and main results}
\label{s_setting_and_main_results}
We consider a lattice system of unbounded continuous spins on the sublattice~$\Lambda = [N] = \{1, \cdots N\} \subset \mathbb{Z}$. The formal Hamiltonian $H:\mathbb{R}^{N} \to \mathbb{R}$ of the system is defined as
\begin{align}\label{e_d_hamiltonian}
H(x) = \sum_{i =1 }^{N} \left( \psi (x_i) + s_i x_i +\frac{1}{2}\sum_{j : \ 1 \leq |j-i| \leq R } M_{ij}x_i x_j \right),
\end{align}
where~$\psi (z ) = \frac{1}{2} z^2 + \psi_b (z)$ and for all $j \notin [N]$ we set $x_j=0$. For each~$i \in [N]$ we define~$M_{ii} : = 1$ and assume the following:
\begin{itemize}
\item The function~$\psi_b: \mathbb{R} \to \mathbb{R}$ satisfies
 \begin{align}\label{e_nonconvexity_bounds_on_perturbation}
 |\psi_b|_{\infty} + |\psi'_b|_{\infty}  + |\psi''_b|_{\infty} < \infty. 
 \end{align}
It is best to imagine~$\psi$ as a double-well potential (see Figure~\ref{f_double_well}).  
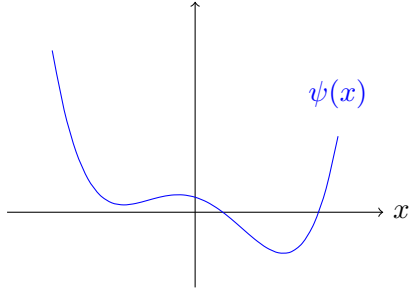
\begin{figure}[t]
\centering
\begin{tikzpicture}
      \draw[->] (-2.5,0) -- (2.5,0) node[right] {$x$};
      \draw[->] (0,-1) -- (0,2.8) node[above] {};
      \draw[scale=1,domain=-1.9:1.9,smooth,variable=\x,blue] plot ({\x},{.3*\x*\x*\x*\x-.7*\x*\x-.3*\x+.2});
      
\node[align=center, below, blue] at (1.9, 1.9) {$\psi(x)$};      
\end{tikzpicture}
\caption{Example of a single-site potential~$\psi$}\label{f_double_well}
\end{figure}

\item The interaction is symmetric, i.e. 
\begin{align}
M_{ij} = M_{ji}  \qquad &\text{for all distinct} \ i, j \in [N]. 
\end{align}

\item The fixed, finite number~$R \in \mathbb{N}$ models the range of interactions between the particles in the system i.e.~it holds that~$M_{ij}=0$ for all~$i,j$ such that~$|i-j|> R$. \\

\item The matrix~$M=(M_{ij})$ is strictly diagonal dominant i.e. for some~$\delta>0$, it holds for any~$i \in [N]$ that
\begin{align} \label{e_strictly_diagonal_dominant}
\sum_{ 1 \leq |j-i| \leq R} |M_{ij}| + \delta \leq M_{ii} = 1.
\end{align}

\item The vector~$s = (s_i) \in \mathbb{R}^{N}$ is arbitrary. It models the interaction with an inhomogeneous external field. Because the interaction is quadratic, this term also models the interaction of the system with boundary values. 
\end{itemize}

% \begin{definition}
% Let~$\Lambda \subset \mathbb{Z}$ be a finite sublattice. We call the spin values~$x^{\mathbb{Z} \setminus \Lambda}$ tempered if
% \begin{align} \label{d_tempered_state}
% \sum_{j \in \mathbb{Z} \setminus \Lambda} |M_{ij}||x_j| < \infty.
% \end{align} 
% \end{definition}

% \begin{remark}
% As we assume the finite range interaction,~\eqref{d_tempered_state} is automatically satisfied for any finite sublattice~$\Lambda \subset \mathbb{Z}$ and spin values~$x^{\mathbb{Z} \setminus \Lambda}$.
% \end{remark}

\begin{definition}

The gce~$\mu^{\sigma}$ associated to the Hamiltonian~$H$ is the probability measure on~$\mathbb{R}^{N}$ given by the Lebesgue density
\begin{align} \label{d_gc_ensemble}
\mu^{ \sigma} \left(dx\right) : = \frac{1}{Z} \exp\left( \sigma \sum_{k =1}^N x_k  - H(x) \right) dx, 
\end{align}
where~$dx$ denotes the Lebesgue measure on~$\mathbb{R}^{N}$. The ce emerges from the gce by conditioning on the mean spin
\begin{align} \label{e_spin_restriction}
  \frac{1}{N} \sum_{k =1}^{N} x_k = m.
\end{align}
More precisely, the ce~$\mu_m $ is the probability measure on
\begin{align}\label{e_d_X_lambda_m}
X_{N, m} : = \left \{ x \in \mathbb{R}^{N} : \ \frac{1}{N} \sum_{k=1}^{N} x_k =m  \right \} \subset \mathbb{R}^{N}
\end{align}
with density
\begin{align} 
\mu_m  (dx) : &= \mu^{\sigma}  \left(dx  \mid \frac{1}{N}\sum_{k =1}^{N} x_k =m \right)  \\
&= \frac{1}{Z} \mathds{1}_{ \left\{ \frac{1}{N} \sum_{k =1}^{N} x_k =m \right\}}\left(x\right) \exp\left( - H(x) \right) \mathcal{L}^{N-1}(dx), \label{d_ce}
\end{align}
where~$\mathcal{L}^{N-1}(dx)$ denotes the~$(N-1)$-dimensional Hausdorff measure supported on~$X_{N, m}$. 
\end{definition}
To relate the external field~$\sigma$ of~$\mu^{\sigma}$ and the mean spin~$m$ of~$\mu_m$ we make the following definition which will be justified in Section~\ref{s_auxiliary_lemmas}.
\begin{definition}
For each~$m \in \mathbb{R}$, we choose~$\sigma = \sigma(m) \in \mathbb{R}$ such that
\begin{align}
\frac{1}{N} \sum_{k=1}^N m_k   = m,
\end{align}
where~$m_k : = \int x_k \mu^{\sigma} (dx)$.
\end{definition}

\begin{definition}[Logarithmic Sobolev inequality (LSI)]
Let~$X$ be a Euclidean space. A Borel probability measure~$\mu$ on~$X$ satisfies the LSI with constant~$\varrho>0$ (or LSI$(\varrho)$) if, for all nonnegative locally Lipschitz functions $f \in L^1 (\mu)$,
\begin{align}\label{e_d_LSI}
\int f \ln f d\mu - \int f d\mu \ln \left( \int f d\mu \right) \leq \frac{1}{2 \varrho} \int \frac{ \left| \nabla f \right| ^2  }{f} d\mu,
\end{align}
where~$\nabla$ denotes the gradient in the Euclidean space~$X$.
\end{definition}
It is well known that on the one-dimensional lattice, the gce satisfies a uniform LSI if the interaction decays fast enough.

\begin{theorem}[Theorem 1.6 in~\cite{MeNi14}]\label{p_uniform_LSI_gce} Let~$H : \mathbb{R}^{N} \to \mathbb{R}$ be the formal Hamiltonian defined as in~\eqref{e_d_hamiltonian}. Assume the interaction is symmetric and the interaction range is infinite. That is,~$M_{ij}$ is not necessarily~$0$ when~$|i-j|>R$. Assume further that the matrix~$M=(M_{ij})$ is strictly diagonal dominant in the sense that there is a~$\delta>0$ with
\begin{align}
\sum_{j : j \neq i} |M_{ij}| + \delta \leq M_{ii}  =1 \qquad \text{for all } i\in [N].
\end{align}
If the interaction decays fast enough, i.e. there are positive constants~$C$ and~$\alpha$ such that
\begin{align}
|M_{ij}| \leq C \frac{1}{|i-j|^{2+\alpha} +1} \qquad \text{for all } i, j \in [N],
\end{align}
then the gce~$\mu^{\sigma}$ satisfies a LSI$(\varrho)$, where~$\varrho>0$ is independent of the system size~$N$ and the external fields~$\sigma$,~$s$. 
\end{theorem}

The main result of this article is that the ce also satisfies on the one-dimensional lattice a uniform LSI for arbitrary strong, finite-range interaction. 
\begin{theorem} \label{p_uniform_lsi}
The ce~$\mu_m$ given by~\eqref{d_ce} satisfies a uniform LSI$(\varrho)$, where~$\varrho>0$ is independent of the system size~$N$, the external field~$s$ and the mean spin~$m \in \mathbb{R}$. 
\end{theorem}

\begin{remark}[From Glauber to Kawasaki] We want to point out that the ce~$\mu_{m}$ is defined on the space~$X_{N,m}$ given by~\eqref{e_d_X_lambda_m}. Because the space~$X_{N,m}$ is endowed with the standard Euclidean structure inherited from~$\mathbb{R}^{N}$, the bound on the right-hand side of the LSI (see~\eqref{e_d_LSI}) is given in terms of the Glauber dynamics. By the discrete Poincar\'e inequality one can recover the bound for the Kawasaki dynamics (cf.~\cite{Cap03} or~\cite[Remark 15]{GrOtViWe09}) in the sense that one endows~$X_{N,m}$ with the Euclidean structure coming from the discrete~$H^{-1}$-norm. The so obtained diffusive scaling of LSI constant for the Kawasaki dynamics is known to be optimal (see~\cite{Yau96} and also Remark 2 in~\cite{Me11}).
\end{remark}
We give the proof of Theorem~\ref{p_uniform_lsi} in Section~\ref{s_proof_uniform_lsi}. There are several basic criteria for the LSI (cf. Appendix~\ref{a_lsi_criteria}), but none of them applies to the ce. The Tensorization Principle~\cite{Gro75} for LSI does not apply because of the restriction to the hyperplane~$X_{N,m}$ and the presence of interaction i.e.~$M\neq0$. The criterion of Bakry \& \'Emery~\cite{BaEm85} does not directly apply because the single-site potentials~$\psi_i$ are non-convex. A combination of the Bakry \& \'Emery and the Holley \& Stroock~\cite{HolStr87} would only yield a LSI with constant~$\varrho$ that is exponentially small in the system size~$N$. The criterion of Otto \& Reznikoff~\cite{OttRez07} does not apply because of the restriction to the hyper-plane~$X_{N,m}$. \\

More advanced methods are needed for deducing a uniform LSI for the ce~$\mu_m$. The most common approach to LSI for Kawasaki dynamics is the Lu \& Yau martingale method~\cite{LuYa93,LaPaYa02,Cha03}. Using this method Landim, Panizo \& Yau~\cite{LaPaYa02} proved Theorem~\ref{p_uniform_lsi} in the special case M = 0 for the Kawasaki bound. An adaptation of this approach by Chafai~\cite{Cha03} lead to the stronger bound for Glauber Dynamics. Providing a new technique, called the two-scale approach, Grunewald, Otto, Westdickenberg (former Reznikoff) \& Villani~\cite{GrOtViWe09} reproduced Theorem~\ref{p_uniform_lsi} for M = 0. In~\cite{Me11}, the uniform LSI was deduced for weak interaction, i.e.~$\|M \| \ll 1$, by perturbing the two-scale approach.\\

The drawback of the two-scale approach is that it elementarily takes advantage of having no interaction term in the Hamiltonian i.e.~setting $M=0$. The basic idea in the two-scale approach is to decompose the lattice~$[N]$ into blocks. This yields a decomposition of the ce~$\mu_{m}$ into a conditional measure, conditioned on the mean spins of the blocks, and a marginal measure for the mean spins. The task is then to deduce two LSIs: A microscopic LSI for the conditional measure and a macroscopic LSI for the marginal measure. After that, the two LSIs are combined into a single LSI for the ce~$\mu_{m}$. If there is no interaction term, blocks do not interact for the conditional and marginal measure. This helps a lot when deducing the microscopic and the macroscopic LSI. If there is a small interaction term, i.e.~$\|M \|\ll 1$, blocks only interact weakly. In~\cite{Me11}, one took advantage of this observation by essentially perturbing the two-scale approach. If there is a large interaction term in the Hamiltonian~$H$ then blocks also interact strongly. It becomes very difficult to deduce the microscopic and macroscopic LSI in the original setting of the two-scale approach. \\

In this article, we overcome this difficulty by using the following strategy. In~\cite{Zeg96}, Zegarlinski deduced the uniform LSI for the gce~$\mu^\sigma$ with strong finite range interaction on the one-dimensional lattice. We follow his approach and decompose the lattice~$[N]$ into odd blocks~$\Lambda_1$ and even blocks~$\Lambda_2$ (see Figure~\ref{f_two_scale_decomposition_of_lattice}). The difference to the two-scale approach is that one does not condition on the mean-spins of blocks but on the spin-values of every even block. The resulting conditional measure~$\mu_m \left(dx^{\Lambda_1} \big| x^{\Lambda_2}\right)$ is also a canonical ensemble but with the advantage that spins in distinct odd blocks do not interact within the Hamiltonian of~$\mu_m \left(dx^{\Lambda_1} \big| x^{\Lambda_2}\right)$ due to the assumption of finite-range interaction. The next step in the Zegarlinski method is to deduce a uniform LSI for~$\mu_m \left(dx^{\Lambda_1} \big| x^{\Lambda_2}\right)$. In our situation this is achieved via the two-scale approach described above. The main new ingredient is a recent result of the authors i.e.~the local Cram\'er theorem (cf.~\cite[Theorem 3]{KwMe17}). The last step in the Zegarlinski method~\cite{Zeg96} is to iteratively condition on the spin values in~$\Lambda_1$ and~$\Lambda_2$ deducing a LSI via a convergence argument. This is where we deviate from the Zegarlinski method. Instead of using an iterative argument we apply the two-scale criterion for LSI (cf.~\cite{OttRez07} or Theorem~\ref{a_two_scale} in the appendix), which in the opinion of the authors is a more direct argument.\\

In the two-scale criterion for the LSI one needs two ingredients: a uniform LSI for the conditional measure~$\mu_m \left(dx^{\Lambda_1} \big| x^{\Lambda_2}\right)$ and a uniform LSI for the marginal measure~$\bar{\mu}_m \left( dx^{\Lambda_2} \right)$. Then, the criterion combines both LSIs into a LSI for the full measure~$\mu_m$. Let us explain how we deduce the LSI for the marginal measure~$\bar{\mu}_m \left( dx^{\Lambda_2} \right)$ which is active on~$\Lambda_2$ and integrates out~$\Lambda_1$. For this, we use the Otto-Reznikoff criterion for LSI (cf.~Theorem~1 in~\cite{OttRez07} or Theorem~\ref{a_otto_reznikoff} in the appendix). The main observation needed is that blocks in~$\Lambda_2$ only interact weakly, if the block size of~$\Lambda_1$ is large enough. We deduce this subtle fact by a series of calculations, building up on the equivalence of ensembles, the decay of correlations, and the moment estimates for the gce. For details of the argument we refer the reader to Section~\ref{s_proof_uniform_lsi}.\\

\section{Auxiliary Lemmas} \label{s_auxiliary_lemmas}

In this section we provide several auxiliary results. All those results were proved in~\cite{KwMe17} for lattice systems with attractive, nearest-neighbor interaction. It is not hard to see that the arguments in~\cite{KwMe17} can be generalized in a straight-forward manner to lattice systems with finite range interaction $R< \infty$, which is considered in this article. The main reason that this is possible is the following. The main technical method used in~\cite{KwMe17} is to  decompose the gce (see definition~\eqref{d_gc_ensemble})
\begin{align}
\mu^{\sigma}(dx) : = \frac{1}{Z} \exp\left( \sigma \sum_{k =1}^{N} x_k -H(x) \right) dx
\end{align}
by conditioning onto even spins~$x_{\tt{even}}$ into
\begin{align}
    \mu^{\sigma}(dx)= \mu^{\sigma}(x_{\tt{odd}} | x_{\tt{even}}) \bar \mu^{\sigma}(d x_{\tt{even}}).
\end{align}
Then one takes advantage that under nearest-neighbor interaction the conditional measure $\mu^{\sigma}(x_{\tt{odd}} | x_{\tt{even}})$ is a product measure. To extend the results of~\cite{KwMe17} to finite-range interaction, one would have to decompose the lattice~$[N]$ into even and odd blocks of size at least the interaction range~$R$. Then instead of conditioning onto the even spins, one would have to condition onto the spins in the even blocks. By the finite range of the interaction the resulting conditional measure would again be a product measure, which is one of the main properties used in~\cite{KwMe17}.  \\

In this section we also show that the results of~\cite{KwMe17} can be extended to arbitrary and not necessarily attractive interaction. More precisely, the interaction~$M_{ij}$ can take on any sign and not only~$M_{ij}\leq 0$. 
In~\cite{KwMe17}, the assumption of attractive interaction, i.e.~$M_{ij}\leq 0$, is only needed to show the lower bound (see Lemma 3 in~\cite{KwMe17})
\begin{align}
\frac{1}{C} \leq \frac{1}{N} \var_{\mu^{ \sigma}} \left( \sum_{k=1}^N X_k \right).
\end{align}

In the next lemma, we provide this lower bound for general and not necessarily attractive interaction. This extends all results of~\cite{KwMe17} to arbitrary interaction. 

\begin{lemma}[Extension of Lemma 3 in~\cite{KwMe17}] \label{l_variance_estimate}
There exists a constant~$C \in (0, \infty)$, uniform in $N$, $s$, and~$\sigma$ such that
\begin{align}\label{e_variance_estimate_mgce}
\frac{1}{C} \leq \frac{1}{N} \var_{\mu^{ \sigma}} \left( \sum_{k=1}^N X_k \right) \leq C.
\end{align}
\end{lemma}

We give the proof of Lemma~\ref{l_variance_estimate} in Section~\ref{s_proof_lemma}. The free energy~$A : \mathbb{R} \to \mathbb{R}$ of the gce~$\mu^{ \sigma}$ is defined by
\begin{align}
A (\sigma) : = \frac{1}{N} \ln \int_{\mathbb{R}^N} \exp \left( \sigma \sum_{k =1}^{N} x_k -H(x) \right) dx.
\end{align}

By Lemma~\ref{l_variance_estimate} we are able to apply Lemma 1 in~\cite{KwMe17} which yields that the free energy~$A$ is uniformly strictly convex. More precisely, it holds

\begin{lemma} [Lemma 1 in~\cite{KwMe17}] \label{l_uniformly_strictly_convex_free_energy}
There exists a constant~$C \in (0, \infty)$ such that for all~$\sigma \in \mathbb{R}$,
\begin{align}
\frac{1}{C} \leq \frac{d^2}{d \sigma^2} A (\sigma) \leq C.
\end{align}
\end{lemma}

With the help of Lemma~\ref{l_uniformly_strictly_convex_free_energy}, we relate the external field~$\sigma$ of~$\mu^{ \sigma}$ and the mean spin~$m$ of~$\mu_m$ as follows:
\begin{definition} \label{a_relation_sigma_m} 
We choose~$\sigma = \sigma(m) \in \mathbb{R}$ and~$m= m(\sigma) \in \mathbb{R}$ such that
\begin{align}
\frac{d}{d\sigma} A(\sigma) = m. \label{e_relation_sigma_m}
\end{align}
Setting~$m_k := \int x_k \mu^{ \sigma}\left(dx\right)$ we equivalently get
\begin{align}
\frac{1}{N} \sum_{k=1}^{N} m_k = \frac{1}{N} \int \left( \sum_{k =1}^N x_k \right) \mu^{ \sigma}\left(dx\right)  = m. \label{e_m=sum_mi}
\end{align}
By strict convexity of~$A(\sigma)$, for each~$m \in \mathbb{R}$ there exists a unique~$\sigma= \sigma(m)$ satisfying~\eqref{e_relation_sigma_m} or vice versa.
\end{definition}

We need the following moment estimates for the gce~$\mu^{\sigma}$.

\begin{lemma}[(23) in~\cite{KwMe17}] \label{l_moment_estimate}
For each~$k \in \mathbb{N}$, there is a constant~$C=C(k)$ such that for each~$i \in [N]$
\begin{align}
\mathbb{E}_{\mu^{ \sigma}}\left[ \left| X_i -m_i \right|^k \right] \leq C(k).
\end{align}
\end{lemma}

\begin{lemma}[Lemma 5 in~\cite{KwMe17}]  \label{l_amgm}
For any finite set~$B_i \subset [N]$ and~$k \in \mathbb{N}$, it holds that
\begin{align}
\left|\mathbb{E}_{\mu^{ \sigma}}\left[ \sum_{i_1 \in B_1 } \cdots \sum_{i_k \in B_k} \left( X_{i_1} - m_{i_1} \right) \cdots \left( X_{i_k} - m_{i_k} \right)  \right]  \right| \lesssim \left| B_1 \right| \cdots \left| B_k \right|.
\end{align}
\end{lemma}

Next, it holds that on the one-dimensional lattice with finite-range interaction,~$\mu^{ \sigma}$ has uniform exponential decay of correlations.

\begin{lemma} [Lemma 6 in~\cite{KwMe17}] \label{l_exponential_decay_of_correlations}
For a function~$f : \mathbb{R}^{N} \to \mathbb{C}$, denote~$\supp f = \{i_1, \cdots, i_k \}$ by the minimal subset of~$[N]$ with~$f(x) = f(x_{i_1}, \cdots, x_{i_k})$. For any functions~$f, g : \mathbb{R}^{N} \to \mathbb{C}$, it holds that
\begin{align}
&\left| \cov_{\mu^{\sigma}} \left( f(X), g(X) \right)\right| \leq C \|\nabla f \|_{L^2 ( \mu^{ \sigma})}\|\nabla g \|_{L^2 ( \mu^{ \sigma})} \exp \left( -C \text{dist} \left( \supp f , \supp g \right) \right).
\end{align}
\end{lemma}

Let~$g$ be the density of the random variable
\begin{align}
\frac{1}{\sqrt{N}} \sum_{k =1}^{N} \left( X_k - m \right) \overset{\eqref{e_m=sum_mi}}{=} \frac{1}{\sqrt{N}} \sum_{k =1}^{N} \left( X_k - m_k \right),
\end{align}
where the random vector~$X=(X_i)_{i \in [N]}$ is distributed according to~$\mu^{ \sigma}$.
The following proposition provides estimates for~$g(0)$.

\begin{proposition} [Proposition 1 in~\cite{KwMe17}] \label{p_main_computation}
For each~$\alpha>0$ and~$\beta > \frac{1}{2}$, there exist uniform constant~$C \in (0, \infty)$ and~$N_0 \in \mathbb{N}$ such that for all~$N \geq N_0$, it holds for all~$\sigma \in \mathbb{R}$ that
\begin{align}
\frac{1}{C} \leq g(0) \leq C , \qquad 
\left| \frac{d}{d\sigma} g(0) \right| \lesssim N^{\alpha} \qquad \text{and} \qquad
\left| \frac{d^2}{d\sigma^2} g(0) \right| \lesssim N^{\beta}.
\end{align}
\end{proposition}
We also need the following statement.
\begin{proposition}[Corollary 2 in~\cite{KwMe17}]\label{p_strict_convexity_bar_h}
Let~$H : \mathbb{R}^N \to \mathbb{R}$ be a Hamiltonian that satisfies the assumptions made in Section~\ref{s_setting_and_main_results}. The coarse-grained Hamiltonian~$\bar H : \mathbb{R} \to \mathbb{R}$ is defined by
\begin{align}
\bar H (m) := - \frac{1}{N} \ln \int_{ \left\{ x \in \mathbb{R}^N \ | \ \frac{1}{N} \sum_{i=1}^N x_i = m \right\} } \exp (-H(x) ) \mathcal{L} (dx). 
\end{align}
Then there is an positive integer~$N_0$ such that for all~$N \geq N_0$ the coarse-grained Hamiltonian~$\bar H$ is uniformly strictly convex. More precisely, there is a uniform constant~$0<C<\infty$ such that for all~$m \in \mathbb{R}$
\begin{align}
  \frac{1}{C} \leq \frac{d^2}{dm^2} \bar H(m) \leq C.
\end{align}
\end{proposition}
The last statement is the decay of correlations of the canonical ensemble.

\begin{proposition}[Theorem 4 in~\cite{KwMe18a}] \label{l_decay_corr_ce}
For each~$f, g : \mathbb{R}^{N} \to \mathbb{R}$, denote~$C(f,g)$ by
\begin{align}
C (f,g) : &= \| \nabla \left( \left(f(X) - \mathbb{E}_{\mu^{\sigma}}\left[ f(X)\right]\right)\left(g(X) - \mathbb{E}_{\mu^{\sigma}}\left[ g(X)\right]\right)   \right)\|_{L^{4}(\mu^{ \sigma})} \\
& \quad + \| \nabla f \|_{L^{4}( \mu^{\sigma})}\| \nabla g \|_{L^{4}(\mu^{ \sigma})}.    \label{e_def_c}
\end{align}
Then for each~$\varepsilon>0$, there exist constants~$N_0$ and~$\tilde{C}=\tilde{C}(\varepsilon) \in (0, \infty)$ independent of the external field~$s = (s_i)$, and the mean spin~$m$ such that for all~$N \geq N_0$, it holds that
\begin{align}
\left|  \cov_{\mu_m}  (f,g) \right| &\leq \tilde{C} \ C(f,g) \left( \frac{  |\supp f | + |\supp g | }{N^{1- \varepsilon}} + \exp\left(-C\text{dist}\left( \supp f, \supp g \right) \right) \right).
\end{align}
\end{proposition}

\subsection{Proof of Lemma~\ref{l_variance_estimate}} \label{s_proof_lemma}

We begin with providing an auxiliary lemma. It is an estimate of single-site variance.

\begin{lemma}\label{p_single_site_variance_bound}
Let~$X\in \mathbb{R}^N$ be a random vector distributed according to~$\mu^{\sigma}$. Then there is a universal constant~$0< C < \infty$ (depending only on the interaction matrix~$M$ and the nonconvexity~$\psi_b$) such that for all~$i \in [N]$ 
\begin{align}
  \label{e_single_site_variance_bound}
  \var_{\mu^{\sigma}} (X_i) \geq C.
\end{align}
\end{lemma}

\noindent \emph{Proof of Lemma~\ref{p_single_site_variance_bound}.} \ By conditioning it hols that
\begin{align}
  \var_{\mu^\sigma} (X_i) & = \mathbb{E}_{\mu^\sigma} \left[ \var_{\mu^\sigma} \left( X_i| (X_j)_{ j \neq i} \right) \right] + \var_{\mu^\sigma} \left( \mathbb{E} \left[ X_i| (X_j)_{ j \neq i} \right]  \right) \\
& \geq \mathbb{E}_{\mu^\sigma} \left[ \var_{\mu^\sigma} \left( X_i| (X_j)_{ j \neq i} \right) \right].
\end{align}
The desired estimate~\eqref{e_single_site_variance_bound} will follow from the uniform bound
\begin{align}
  \var_{\mu^\sigma} \left( X_i| (X_j)_{ j \neq i} \right)  \geq C>0.
\end{align}
Indeed, the conditional measure~$\mu^\sigma \left( dx_i| (x_j)_{ j \neq i} \right) $ has the Lebesgue density
\begin{align}
 \mu^\sigma \left( dx_i| (x_j)_{ j \neq i} \right)= \frac{1}{Z} \exp \left( - \frac{1}{2} x_i^2 -\tilde s_i x_i - \psi_b(x_i) \right),
\end{align}
where~$\tilde s_i$ is given by
\begin{align}
  \tilde s_i = s_i + \frac{1}{2} \sum_{j \neq i} M_{ij} x_j. 
\end{align}
Let~$\nu$ denote the one-dimensional measure given by the Lebesgue density
\begin{align}
  \nu (dz) = \frac{1}{Z_\nu} \exp \left( - \frac{1}{2} z^2 - \tilde s z \right) dz
\end{align}
Using the bound~$|\psi_b|_\infty\leq C$ and the optimality of the mean for the variance, and the fact that~$\var_{\nu}(Z)=1$ we obtain the desired estimate
\begin{align}
  \var_{\mu^\sigma} \left( X_i| (X_j)_{ j \neq i} \right)  & \geq  \exp \left(-2C \right)  \frac{1}{Z_\nu} \int \left( z - \mathbb{E}_{\nu} \left[X_i  \right] \right)^2  \exp \left(- \frac{1}{2} z^2 - \tilde s z \right) dz \\
& \geq \exp \left(-2C \right) \var_{\nu}(Z) = \exp \left(-2C \right) =C>0.
\end{align}
\qed
\medskip

Let us now see how the Lemma~\ref{p_single_site_variance_bound} yields Lemma~\ref{l_variance_estimate}. \\

\noindent \emph{Proof of Lemma~\ref{l_variance_estimate}.} \ For the proof of the upper bound, we refer to the proof of~\cite[Lemma 3]{KwMe17}. Here, we shall only provide the proof of the lower bound in~\eqref{e_variance_estimate_mgce}. Let~$Q$ be the set~$\left\{ R+1, 2(R+1), 3(R+1), \cdots \right\} \cap [N]$. By conditioning we get that
\begin{align}
 \var_{\mu^\sigma} \left( \sum_{i \in [N] }X_i \right) & = \mathbb{E}_{\mu^\sigma} \left[ \var_{\mu^\sigma} \left( \sum_{i \in Q } X_i \bigg| (X_j)_{ j \in [N] \backslash Q} \right) \right] \\
& \quad + \var_{\mu^\sigma} \left( \mathbb{E} \left[ \sum_{i \in Q} X_i  \bigg| (X_j)_{ j \in [N] \backslash Q } \right]  \right) \\
& \geq \mathbb{E}_{\mu^\sigma} \left[ \var_{\mu^\sigma} \left( \sum_{i \in Q} X_i \bigg| (X_j)_{ j \in [N] \backslash Q} \right) \right] \\
& = \sum_{i \in Q}\mathbb{E}_{\mu^\sigma} \left[ \var_{\mu^\sigma} \left( X_i \bigg| (X_j)_{ j \in [N] \backslash Q} \right) \right] ,
\end{align}
where we used in the last line the fact that because the interaction range is~$R$, different sites in~$Q$ become independent after conditioning onto the spin values~$(X_j)_{ j \in [N] \backslash Q} $. Now, an application of Lemma~\ref{p_single_site_variance_bound} yields that
\begin{align}
  \var_{\mu^\sigma} \left( X_i \right) & \geq  \sum_{i \in Q}\mathbb{E}_{\mu^\sigma} \left[ \var_{\mu^\sigma} \left( X_i | (X_j)_{ j \in [N] \backslash Q} \right) \right] \\
&\geq   C \ |Q| = \frac{C}{R+1} N,
\end{align}
which is the desired statement.
\qed

\medskip

\medskip

\section{Proof of Theorem~\ref{p_uniform_lsi}} \label{s_proof_uniform_lsi}
In the proof of Theorem~\ref{p_uniform_lsi} we use an idea of Zegarlinski and decompose the lattice~$[N]$ into two parts. This idea was used to prove the uniform LSI for the gce~$\mu^\sigma$ on the one-dimensional lattice (see~\cite{Zeg96}). However, instead of using an iterative argument as in~\cite{Zeg96} we make use of ideas outlined in~\cite{OttRez07} and~\cite{GrOtViWe09}.\\

We decompose the lattice into two types of blocks (see Figure~\ref{f_two_scale_decomposition_of_lattice}):
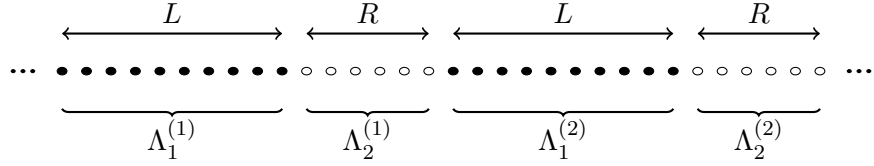
\begin{figure}[t]
\centering
\begin{tikzpicture}[xscale=1.3]

\draw[fill] (0,0) circle [radius=0.02];
\draw[fill] (.1,0) circle [radius=0.02];
\draw[fill] (.2,0) circle [radius=0.02];

\draw[fill] (.5,0) circle [radius=0.05];
\draw[fill] (.75,0) circle [radius=0.05];
\draw[fill] (1,0) circle [radius=0.05];
\draw[fill] (1.25,0) circle [radius=0.05];
\draw[fill] (1.5,0) circle [radius=0.05];
\draw[fill] (1.75,0) circle [radius=0.05];
\draw[fill] (2,0) circle[radius=0.05];
\draw[fill] (2.25,0) circle[radius=0.05];
\draw[fill] (2.5,0) circle[radius=0.05];
\draw[fill] (2.75,0) circle[radius=0.05];

\draw (3,0) circle[radius=0.05];
\draw (3.25,0) circle[radius=0.05];
\draw (3.5,0) circle[radius=0.05];
\draw (3.75,0) circle[radius=0.05];
\draw (4,0) circle[radius=0.05];
\draw (4.25,0) circle[radius=0.05];

\draw[fill] (4.5,0) circle [radius=0.05];
\draw[fill] (4.75,0) circle [radius=0.05];
\draw[fill] (5,0) circle [radius=0.05];
\draw[fill] (5.25,0) circle [radius=0.05];
\draw[fill] (5.5,0) circle [radius=0.05];
\draw[fill] (5.75,0) circle [radius=0.05];
\draw[fill] (6,0) circle[radius=0.05];
\draw[fill] (6.25,0) circle[radius=0.05];
\draw[fill] (6.5,0) circle[radius=0.05];
\draw[fill] (6.75,0) circle[radius=0.05];

\draw (7,0) circle[radius=0.05];
\draw (7.25,0) circle[radius=0.05];
\draw (7.5,0) circle[radius=0.05];
\draw (7.75,0) circle[radius=0.05];
\draw (8,0) circle[radius=0.05];
\draw (8.25,0) circle[radius=0.05];

\draw[fill] (8.55,0) circle [radius=0.02];
\draw[fill] (8.65,0) circle [radius=0.02];
\draw[fill] (8.75,0) circle [radius=0.02];

\draw[decorate,decoration={brace,mirror},thick] (0.5,-.5) -- node[below]{$\Lambda_1 ^{(1)}$} (2.75,-.5);
\draw[decorate,decoration={brace,mirror},thick] (3,-.5) -- node[below]{$\Lambda_2 ^{(1)}$} (4.25,-.5);
\draw[decorate,decoration={brace,mirror},thick] (4.5,-.5) -- node[below]{$\Lambda_1 ^{(2)}$} (6.75,-.5);
\draw[decorate,decoration={brace,mirror},thick] (7,-.5) -- node[below]{$\Lambda_2^{(2)}$} (8.25,-.5);

\draw[thick,<->] (.5,.5) -- (2.75,.5);
\node[align=center, above] at (1.625,.5) {$L$};
\draw[thick,<->] (3,.5) -- (4.25,.5);
\node[align=center, above] at (3.625,.5) {$R$};
\draw[thick,<->] (4.5,.5) -- (6.75,.5);
\node[align=center, above] at (5.625,.5) {$L$};
\draw[thick,<->] (7,.5) -- (8.25,.5);
\node[align=center, above] at (7.625,.5) {$R$};

\end{tikzpicture}
\caption{Arrangement in the cell~$[1, 2(L+R)]$ for~$L=10$ and~$R=6$ }\label{f_two_scale_decomposition_of_lattice}
\end{figure}

\begin{align}
\Lambda &: = [N] = \{1, 2, \cdots, N \}, \\
\Lambda_{1} &: = \bigcup_{ n \in \mathbb{Z}} \Lambda \cap \left( \left[1, L \right]  + \left(L+R\right)(n-1) \right) = \bigcup_{n =1}^{M} \Lambda_1^{(n)} \label{d_lambda1}, \\
 \Lambda_{2} &: = \bigcup_{ n \in \mathbb{Z}} \Lambda \cap  \left( \left[ L+1 , L+R \right] + \left(L+R \right)(n-1) \right) = \bigcup_{n =1}^{M} \Lambda_2 ^{(n)}, \label{d_lambda2}
\end{align}
where~$R$ is the range of interactions between particles (cf.~\eqref{e_d_hamiltonian}). The number~$L$ will be chosen later. \\

Recall the definition~\eqref{d_ce} of the ce~$\mu_m$
\begin{align}
\mu_m (dx) =  \frac{1}{Z} \mathds{1}_{ \left\{ \frac{1}{N} \sum_{k=1}^N x_k =m \right\}}(x)\exp\left(  - H\left(x\right) \right) \mathcal{L}^{N-1}(dx).
\end{align}
The decomposition of~$\Lambda=\Lambda_1 \cup \Lambda_2$ into odd blocks given by~$\Lambda_1$ and even blocks given by~$\Lambda_2$ yields a decomposition of the ce~$\mu_m$ into conditional measure~$\mu_m \left( dx^{\Lambda_1} \big| x^{ \Lambda_2} \right)$ and a marginal measure $\bar{\mu}_m \left( dx^{ \Lambda_2} \right)$. That is, for any test function~$f$, it holds that
\begin{align}\label{e_zegarlinski_decomposition}
\int f(x) \mu_m (dx) = \int \left( \int f(x) \mu_m \left( dx^{\Lambda_1} \big| x^{ \Lambda_2} \right) \right)\bar{\mu}_m \left( dx^{ \Lambda_2} \right).
\end{align}

Now, the strategy is to deduce uniform LSIs for the conditional measure~$\mu_m \left( dx^{\Lambda_1} \big| x^{ \Lambda_2} \right)$ and the marginal measure~$\bar{\mu}_m \left( dx^{ \Lambda_2} \right)$. The uniform LSI for the full measure~$\mu_m$ is then deduced via the two-scale criterion for the LSI (see eg.~\cite{OttRez07} or~\cite{GrOtViWe09}).\\

Let us explain how the uniform LSI for the conditional measures~$\mu_m \left( dx^{\Lambda_1} \big| x^{ \Lambda_2} \right)$ is deduced. By conditioning onto the even blocks~$\Lambda_2$, the spins in one odd block (say~$x_i \in \Lambda_1 ^{(n)}$) do not interact with the spins in another odd block (say~$x_j \in \Lambda_2 ^{(l)}$, $l\neq n$) within the Hamiltonian of~$\mu_m \left( dx^{\Lambda_1} \big| x^{ \Lambda_2} \right)$. The only way they interact is through the constraint~\eqref{e_spin_restriction}. With this observation one can modify the proof of~\cite[Theorem 14]{GrOtViWe09} (see also Theorem~\ref{a_two_scale}) to deduce a uniform LSI for the conditional measure~$\mu_m \left(dx^{\Lambda_1} \big| x^{\Lambda_2} \right)$ via the two-scale approach (see Proposition~\ref{p_conditional_lsi} below). \\

Let us turn to the uniform LSI for the marginal measure~$\bar{\mu}_m$. We observe that the marginal measure~$\bar{\mu}_m$ is not constrained onto a hyperplane. Because~$\bar{\mu}_m$ can be interpreted as a gce on a one-dimensional lattice, the marginal measure~$\bar{\mu}_m$ should heuristically satisfy a uniform LSI. Rigorously, the uniform LSI for the marginal measure~$\bar{\mu}_m$ is deduced via the the Otto-Reznikoff criterion for LSI (cf. Theorem~\ref{a_otto_reznikoff}). For this we need two ingredients: The first ingredient is to show that on each block~$\Lambda_2 ^{(n)}$, the conditional measures~$\bar{\mu}_m \left( dx^{\Lambda_2 ^{(n)}} \big| x^{\Lambda_2 ^{(l)}}, l \neq n \right)$ satisfy a LSI with a constant that is uniform in the conditioned values~$x^{\Lambda_2 ^{(l)}}, l \neq n$. This ingredient is derived via an adaptation of the argument of~\cite{Me13}. The main ingredients for this part are the local Cram\'er theorem (see~\cite[Theorem 3]{KwMe17}), the Holley-Stroock Principle (see Theorem~\ref{a_holley_stroock}), the decay of correlations of the gce (see Lemma~\ref{l_exponential_decay_of_correlations}) and a decomposition method for Hamiltonians introduced in~\cite{Me13}. The second ingredient is to show that, in the Hamiltonian of~$\bar{\mu}_m$, the interactions between blocks~$\Lambda_2 ^{(l)}$ and~$\Lambda_2 ^{(n)}$,~$l \neq n$, are sufficiently small. This is achieved by observing that the correlations of the conditional measure~$\mu_m \left(dx^{\Lambda_1} \big| x^{\Lambda_2} \right)$ decay due to~\cite[Theorem 4]{KwMe18a} and choosing the parameter~$L$ sufficiently large. For more details, see Proposition~\ref{p_marginal_lsi} from below. \\ 

Now, let us turn to the detailed argument. We prove the following three propositions. The first one is the uniform LSI for the conditional measure~$\mu_m \left( dx^{\Lambda_1} \big| x^{ \Lambda_2} \right)$.
\begin{proposition} \label{p_conditional_lsi}
The conditional measure~$\mu_m \left( dx^{\Lambda_1} \big| x^{ \Lambda_2} \right)$ satisfies a LSI with constant~$\varrho_1 >0$ that is uniform in~$m$,~$s$ and the conditioned spins~$x^{ \Lambda_{2}}$. The constant~$\varrho_1 >0$ is uniform in~$K=|\Lambda_1|$,~$x^{\Lambda_2}$,~$s$ and~$m$.
\end{proposition}
The second proposition states the LSI for the marginal measure~$\bar{\mu}_m \left( dx^{ \Lambda_2} \right)$.
\begin{proposition} \label{p_marginal_lsi}
There is a constant~$L_0$ such that if~$L \geq L_0$ the marginal measure~$\bar{\mu}_m \left( dx^{ \Lambda_2} \right)$ satisfies a LSI with constant~$\varrho_2>0$ that is uniform in~$m$,~$s$ and~$L$. The constant~$\varrho_2$ only depends on the interaction range~$R$.
\end{proposition}
The last proposition is
\begin{proposition} \label{p_combined_lsi}
Let~$\mu_m$ be the ce defined by~\eqref{d_ce}. Assume that
\begin{itemize}
\item The conditional measure~$\mu_m \left( dx^{\Lambda_1} \big| x^{ \Lambda_2} \right)$ satisfies LSI$(\varrho_1)$. The constant~$\varrho_1$ is uniform in~$m$,~$s$ and the conditioned spins~$x^{\Lambda_2}$.
\item The marginal measure~$\bar{\mu}_m \left( dx^{ \Lambda_2} \right)$ satisfies LSI$(\varrho_2)$. The constant~$\varrho_2$ is uniform in~$m$, $s$ and~$L$.
\end{itemize}
Then the ce~$\mu_m$ satisfies LSI$(\varrho)$ with a constant~$\varrho>0$ that is uniform in~$N$,~$m$ and~$s$.
\end{proposition}
The main work of this article consists of deducing Proposition~\ref{p_conditional_lsi} and Proposition~\ref{p_marginal_lsi} which are done in Section~\ref{s_proof_conditional_lsi} and Section~\ref{s_proof_marginal_lsi} respectively. Proposition~\ref{p_combined_lsi} is a slight modification of~\cite[Theorem 3]{GrOtViWe09}. For convenience of the readers, we give the proof of Proposition~\ref{p_combined_lsi} in Section~\ref{s_proof_combined_lsi}. \\

\noindent \emph{Proof of Theorem~\ref{p_uniform_lsi}.} \ Without loss of generality, it is sufficient to prove the uniform LSI for large enough system size~$N \geq N_0$. For small system size~$N \leq N_0$ the LSI can be deduced via a combination of Bakry-\'Emery criterion (cf. Theorem~\ref{a_bakry_emery}) and the Holley-Stroock perturbation principle (see Theorem~\ref{a_holley_stroock}). The uniform LSI for large systems~$N \geq N_0$ is obtained by choosing~$L$ large but fixed and a combination of Proposition~\ref{p_conditional_lsi}, Proposition~\ref{p_marginal_lsi} and Proposition~\ref{p_combined_lsi}. 
\qed

\section{Proof of Proposition~\ref{p_conditional_lsi}} \label{s_proof_conditional_lsi}

The goal is to deduce a uniform LSI for the conditional measure~$\mu_m \left(dx^{\Lambda_1} \big| x^{\Lambda_2} \right)$. As described in Section~\ref{s_proof_uniform_lsi}, the proof of Proposition~\ref{p_conditional_lsi} is an adaptation of~\cite[Theorem 14]{GrOtViWe09}, where a uniform LSI was deduced for a ce with non-interacting Hamiltonian. As explained in Section~\ref{s_setting_and_main_results} the main observation is that by conditioning on even blocks, odd blocks do not interact within the Hamiltonian and the setting of~\cite {GrOtViWe09} applies with minor adaptation. More precisely, let us fix~$x^{\Lambda_2} \in \mathbb{R}^{\Lambda_2}$ and denote~$K := |\Lambda_1|$. We can rewrite the spin restriction~\eqref{e_spin_restriction} as
\begin{align}
\frac{1}{K} \sum_{k \in \Lambda_1} x_k = \frac{N}{K} m - \frac{1}{K} \sum_{i \in \Lambda_2} x_i =: \widetilde{m}. \label{e_def_tilde_m}
\end{align}
Assume that we have~$M$ blocks in~$\Lambda_1$. That is,
\begin{align}
\Lambda_1 = \bigcup_{n=1}^{M} \Lambda_1^{(n)}.
\end{align}
Then the Hamiltonian~$H$ can be written as
\begin{align}\label{e_decomposition_Hamiltonian_lambda_1}
H(x) = \sum_{n=1}^M H_n \left(x^{\Lambda_1 ^{(n)}} \right) + C\left(x^{\Lambda_2} \right),
\end{align}
where
\begin{align}
H_n \left(x^{\Lambda_1 ^{(n)}} \right) = \sum_{i \in \Lambda_1 ^{(n)}} \left( \psi (x_i) + s_i x_i + \frac{1}{2} \sum_{j: \ 1 \leq |j-i| \leq R } M_{ij} x_i x_j \right),
\end{align}
and~$C\left( x^{\Lambda_2} \right)$ is a constant that only depends on~$x^{\Lambda_2}$. We observe that for each~$n$ the potential~$H_n \left(x^{\Lambda_1 ^{(n)}} \right)$ only depends on~$x^{\Lambda_1}$ through the spin values~$x^{\Lambda_1 ^{(n)}}$. The potential~$H_n \left(x^{\Lambda_1 ^{(n)}} \right)$ depends on the conditioned spin values~$x^{\Lambda_2}$ only via a linear term. Most importantly, there is no interaction between odd blocks~$\left\{ \Lambda_1^{(n)} : n=1, \cdots, M \right\}$ within the Hamiltonian~$H$. \\

Let us now explain how the argument from~\cite[Theorem 14]{GrOtViWe09} applies. We only point out the main differences and leave the details as an exercise. We start with observing that the conditional measure $\mu_m \left(dx^{\Lambda_1} \big| x^{\Lambda_2} \right)$ is a ce wrt.~the fixed mean spin~$\widetilde{m}$ given by~\eqref{e_def_tilde_m}. Let~$L$ denote the size of a block in~$\Lambda_1$ i.e.~$L = |\Lambda_1^{(1)}|$. Let~$P$ be defined as the map
\begin{align}
  P x^{\Lambda_1} = \left( \frac{1}{L} \sum_{i \in \Lambda_1^{(n)}} x_i \right)_{n \in [M]} =: (y_n)_{n \in [M]} =:y
\end{align}
that associates to every block~$\Lambda_1^{(n)}$ it's mean spin~$y_n$. The mapping~$P$ yields a decomposition of the measure~$\mu_m \left(dx^{\Lambda_1} \big| x^{\Lambda_2} \right)$ into a conditional measure~$\mu_m \left(dx^{\Lambda_1} \big| Px^{\Lambda_1}=y, \  x^{\Lambda_2} \right)$ and a marginal measure~$\bar \mu_m \left(dy \big| \  x^{\Lambda_2} \right)$ i.e.
\begin{align}
  \mu_m \left(dx^{\Lambda_1} \big| x^{\Lambda_2} \right) = \mu_m \left(dx^{\Lambda_1} \big| Px^{\Lambda_1}=y, \  x^{\Lambda_2} \right) \ \bar \mu_m \left(dy \big| \  x^{\Lambda_2} \right).
\end{align}
The core of the argument in~\cite{GrOtViWe09} is to deduce a LSI for the measure $ \mu_m \left(dx^{\Lambda_1} \big| Px^{\Lambda_1}=y, \  x^{\Lambda_2} \right) $ and for~$\bar \mu_m \left(dy \big| \  x^{\Lambda_2} \right)$. Those two LSIs then combine into a LSI for the full measure~$\mu_m \left(dx^{\Lambda_1} \big| x^{\Lambda_2} \right)$.\\

Let us consider the LSI for the conditional measure~$ \mu_m \left(dx^{\Lambda_1} \big| Px^{\Lambda_1}=y, \  x^{\Lambda_2} \right)$. We observe that due to~\eqref{e_decomposition_Hamiltonian_lambda_1} and the conditioning onto~$x^{\Lambda_2}$ the measure~$ \mu_m \left(dx^{\Lambda_1} \big| Px^{\Lambda_1}=y, \  x^{\Lambda_2} \right)$ is a product measure on the blocks. It follows that the measure $ \mu_m \left(dx^{\Lambda_1} \big| Px^{\Lambda_1}=y, \  x^{\Lambda_2} \right)$ satisfies a LSI via a combination of the Bakry-\'Emery Criterion (cf.~Theorem~\ref{a_bakry_emery}), the Holley-Stroock Perturbation Principle (cf.~Theorem~\ref{a_holley_stroock} ) and the Tensorization Principle (cf.~Theorem~\ref{a_tensorization}). Because the conditioning~$x^{\Lambda_2}$ only enters the Hamiltonian of~$ \mu_m \left(dx^{\Lambda_1} \big| Px^{\Lambda_1}=y, \  x^{\Lambda_2} \right)$ as a linear term, the obtained LSI constant is uniform in the conditioned values~$x^{\Lambda_2}$, the mean spin~$m$, the linear term~$s$ and the overall system size~$N$. The constant may depend on the size of the odd blocks. Let us turn to the LSI for the marginal measure~$\bar \mu_m \left(dy \big| \  x^{\Lambda_2} \right)$. We use the same strategy as in the proof of~\cite[Theorem 14]{GrOtViWe09}. We observe that the Hamiltonian~$\bar H(y)$ of $\bar \mu_m \left(dy \big| \  x^{\Lambda_2} \right)$ can be written after cancellation of constant terms as
\begin{align}
  \bar H(y) = L \sum_{n=1}^M \bar \psi_n(y_n),  
\end{align}
where the function~$\bar \psi_n$ is given by
\begin{align}
  \bar \psi_n (z) = - \frac{1}{L} \ln \int_{\left\{ x^{\Lambda_1^{(n)}} \in \mathbb{R}^{\Lambda_1^{(n)}} \ | \ \frac{1}{K} \sum_{i \in \Lambda_1^{(n)}} x_i = z \right\} } \exp (- H_n (x^{\Lambda_1^{(n)}})) \mathcal{L} (dx^{\Lambda_1^{(n)}}).
\end{align}
We observe that~$H_n (x^{\Lambda_1^{(n)}}) $ satisfies the same structural assumptions as the Hamiltonian~$H$ in Section~\ref{s_setting_and_main_results}. Hence, an application of Proposition~\ref{p_strict_convexity_bar_h} yields that the function~$\bar \psi_n$ is uniformly strictly convex for large enough~$L$. Hence, it follows that~$\bar H(y) $ is uniformly strictly convex and therefore the marginal measure $\bar \mu_m \left(dy \big| \  x^{\Lambda_2} \right)$ satisfies a uniform LSI by the Bakry-\'Emery criterion (cf.~Theorem~\ref{a_bakry_emery}). Again, the LSI constant is uniform in the conditioned values~$x^{\Lambda_2}$. It is left to combine both LSIs to a single LSI for the measure~$\mu_m \left(dx^{\Lambda_1} \big| x^{\Lambda_2} \right)$ which is done in a similar fashion as in~\cite{{GrOtViWe09}}.

% {\color{blue} Explain main strategy of the proof of~\cite[Theorem 14]{GrOtViWe09}:
% \begin{itemize}
% \item decomposition into blocks.
% \item Strategy is to use two scale approach.
% \item microscopic LSI due to a combination of Holley-Stroock perturbation and Tensorization
% \item macroscopic LSI. Needs observation that coarse grained Hamiltonian is uniformly striclty convex therefore LSI due to Bary-Emery criterion.
% \end{itemize}

% We can apply a similar argument for the proof of Proposition~\ref{p_conditional_lsi} because:
% \begin{itemize}
% \item Due to the fact that the interaction is of finite range~$R$ blocks in~$\Lambda_1$ do not interact.
% \item Thereofore one can deduce the microscopic LSI via the Tensorization, Bakry-\'Emery principle and Holley-Stroock perturbation.
% \item Second main ingredient is the uniform strict convexity of the coarse-grained Hamiltonian. This fact follows from the results in~\cite{KwMe17} due to the fact that blocks do not interact.
% \item Strict convexity of coarse-grained Hamiltonian of full system....unkown difficult but needed when deducing the hydrodynamic limit of the Kawaski dynamic....open problem.
% \end{itemize}

\section{Proof of Proposition~\ref{p_marginal_lsi} } \label{s_proof_marginal_lsi}

The goal is to deduce a uniform LSI for the marginal measure~$\bar \mu_m (dx^{\Lambda_2})$ given by~\eqref{e_zegarlinski_decomposition}. We recall the decomposition of~$\Lambda = \Lambda_1 \cup \Lambda_2$ into odd blocks given by~$\Lambda_1$ and even blocks given by~$\Lambda_2$ (cf.~Figure~\ref{f_two_scale_decomposition_of_lattice}). The marginal measure~$\bar \mu_m (dx^{\Lambda_2})$ arises from~$\mu_m$ by integrating out the spins located in the odd blocks~$\Lambda_1$. We will deduce the uniform LSI for~$\bar \mu_m (dx^{\Lambda_2})$ by applying the Otto-Reznikoff criterion~\cite{OttRez07} (see Theorem~\ref{a_otto_reznikoff} in the appendix). Applying this criterion needs two ingredients. For explaining the first ingredient let us consider the conditional measures $\bar{\mu}_m \left(dx^{\Lambda_2^{(n)}} \big| x^{\Lambda_2^{(l)}} , l \neq n \right)$. This measure arises from the ce~$\mu_m$ by conditioning on all spins in even blocks except of one and integrating out all spins in odd blocks. The first ingredient is a LSI for the measures~$\bar{\mu}_m \left(dx^{\Lambda_2^{(n)}} \big| x^{\Lambda_2^{(l)}} , l \neq n \right)$ that is uniform in the conditioned spins~$x^{\Lambda_2^{(l)}}$, $l \neq n$ (see Lemma~\ref{l_block_marginal_lsi} below). The second ingredient is that the interactions in the measure~$\bar \mu_m(dx^{\Lambda_2})$ between even blocks are sufficiently small (see Lemma~\ref{l_marginal_otto_reznikoff_argument}).\\

The first ingredient looks innocent on the first sight. By integrating out the odd blocks, the marginal measure~$\bar \mu_m (dx^{\Lambda_2})$ is not restricted to a hyperplane anymore. Therefore a LSI should hold because one considers a gce on a one-dimensional lattice. The difficult part is to show that the LSI constant is uniform in the conditioned spin values~$x^{\Lambda_2^{(l)}}$, $l \neq n$. Those values enter the Hamiltonian of the measure~$\bar \mu_m (dx^{\Lambda_2})$ via a subtle interaction term, namely the free energy of the ce~$\mu_m (dx^{\Lambda_1}|x^{\Lambda_2})$. This interaction term is non-quadratic and has infinite range and therefore is hard to control.\\

We derive the uniform LSI by extending a method of~\cite{Me13} for gces of one variable to ces of multiple variables (see (cf.~\cite[Lemma 3.1, Lemma 3.2]{Me13})). In the first step of the argument we use the Holley-Stroock Perturbation Principle (see Theorem~\ref{a_holley_stroock}) to modify the interaction term in the Hamiltonian. By the equivalence of ensembles (see~\cite[Theorem 2]{KwMe17}) this allows us to consider in the Hamiltonian an interaction term that is the free energy of the associated gce, which is easier to control than the free energy of the ce. Then we follow the approach of~\cite{Me13} and write the Hamiltonian~$H$ on a block~$\Lambda_2^{(n)}$ as a bounded perturbation of a uniform strictly convex potential and deduce the uniform LSI via a combination of the Bakry-\'Emery criterion (see~Theorem~\ref{a_bakry_emery}) and the Holley-Stroock Perturbation Principle (see Theorem~\ref{a_holley_stroock}).\\

In the next two lemmas we formulate the main ingredients for the proof of Proposition~\ref{p_marginal_lsi}. They correspond to~\cite[Lemma 3.1, Lemma 3.2]{Me13} with the difference that the lemmas in this article are more general in the sense that we consider multi-variate measures and that the interaction term is more complicated. 
\begin{lemma}\label{l_block_marginal_lsi}
Assume~$K= |\Lambda_1|$ is large enough. Then for each~$n \in \{1, \cdots, M\}$, the block marginal measure~$\bar{\mu}_m \left(dx^{\Lambda_2^{(n)}} \big| x^{\Lambda_2^{(l)}} , l \neq n \right)$ satisfies the LSI$(\tau)$ for some~$\tau>0$ that depends on~$R$ but is independent of the mean spin~$m$, the external field~$s$ and the conditioned spins~$x^{\Lambda_2^{(l)}}, l \neq n$.
\end{lemma}
\begin{lemma}\label{l_marginal_otto_reznikoff_argument}
Let~$\bar{H}_{\Lambda_2}$ be the Hamiltonian of the marginal measure~$\bar{\mu}_m \left( dx^{\Lambda_2} \right)$, i.e.
\begin{align}
\bar{H}_{\Lambda_2}\left(x^{\Lambda_2 }  \right) &= -\ln \int_{\frac{1}{N} \sum_{k \in \Lambda_1} x_k = m - \frac{1}{N} \sum_{i \in \Lambda_2} x_i} \exp\left( - H\left(x^{\Lambda_1}, x^{\Lambda_2}\right) \right)  \mathcal{L}^{K-1}(dx^{\Lambda_1}) \label{d_def_marginal_hamiltonian1} \\
& \overset{\eqref{e_def_tilde_m}}{=} -\ln \int_{\frac{1}{K} \sum_{k \in \Lambda_1} x_k = \widetilde{m}} \exp\left( - H\left(x^{\Lambda_1}, x^{\Lambda_2}\right) \right)  \mathcal{L}^{K-1}(dx^{\Lambda_1}). \label{d_def_marginal_hamiltonian2}
\end{align}
Let~$\varepsilon>0$ be given. Then for any~$n, l \in \{1, \cdots M\}$ with~$n\neq l$, There is a constant~$\tilde{C}=\tilde{C}(\varepsilon)$ such that
\begin{align}
\left \| \left( \nabla_i \nabla_j \bar{H}_{\Lambda_2} \left( x^{\Lambda_2} \right)  \right)_{i \in \Lambda_2 ^{(n)}, \ j \in \Lambda_2 ^{(l)}} \right\| \leq \tilde{C}\frac{R}{K^{1-\varepsilon}} + \tilde{C}R \exp \left( -CL|n-l| \right) .
\end{align}
Here,~$\| \cdot \|$ denotes the operator norm of a bilinear form and~$L$ denotes the size of a single odd block and~$K$ denotes the overall size of all combined odd blocks, i.e.~$K= |\Lambda_1| =  LM$. 
\end{lemma}
The Lemma~\ref{l_block_marginal_lsi} and Lemma~\ref{l_marginal_otto_reznikoff_argument} are proved in Section~\ref{s_block_marginal} and Section~\ref{s_marginal_otto}, respectively.\\

\noindent \emph{Proof of Proposition~\ref{p_marginal_lsi}.} \ Recall that~$K= LM$. By choosing~$L$ large enough, a combination of Lemma~\ref{l_block_marginal_lsi}, Lemma~\ref{l_marginal_otto_reznikoff_argument} and the Otto-Reznikoff criterion (cf. Theorem~\ref{a_otto_reznikoff}) yields that the marginal measure~$\bar{\mu}_m$ satisfies a LSI with the desired properties.\qed

\subsection{Proof of Lemma~\ref{l_block_marginal_lsi} } \label{s_block_marginal}
We note that the block Hamiltonian~$\bar{H}_{\Lambda_2^{(n)}}$ of the block marginal measure $\bar{\mu}_m \left(dx^{\Lambda_2^{(n)}} \big| x^{\Lambda_2^{(l)}} , l \neq n \right)$ is written as follows:
\begin{align}
\bar{H}_{\Lambda_2^{(n)}}\left(x^{\Lambda_2 ^{(n)}} \big| x^{\Lambda_2 ^{(l)}}, \ l \neq n \right) &= -\ln \int_{\frac{1}{N} \sum_{k \in \Lambda_1} x_k = m - \frac{1}{N} \sum_{i \in \Lambda_2} x_i} \exp\left( - H\left(x^{\Lambda_1}, x^{\Lambda_2}\right) \right)  \mathcal{L}(dx^{\Lambda_1}) \\
& \overset{\eqref{e_def_tilde_m}}{=} -\ln \int_{\frac{1}{K} \sum_{k \in \Lambda_1} x_k = \widetilde{m}} \exp\left( - H\left(x^{\Lambda_1}, x^{\Lambda_2}\right) \right)  \mathcal{L}(dx^{\Lambda_1}).
\end{align}
For notational convenience we abbreviate
\begin{align}
\bar{H}_{\Lambda_2^{(n)}} \left( x^{ \Lambda_2 ^{(n)}} \right) = \bar{H}_{\Lambda_2^{(n)}}\left(x^{\Lambda_2 ^{(n)}} \big| x^{\Lambda_2 ^{(l)}}, \ l \neq n \right).
\end{align}
Our aim is to decompose the block Hamiltonian~$\bar{H}_{\Lambda_2^{(n)}}$ into a strictly convex function and a bounded function. Then the proposition will follow from the Bakry-\'{E}mery criterion (cf. Theorem~\ref{a_bakry_emery}) and the Holley-Stroock Perturbation Principle (cf. Theorem~\ref{a_holley_stroock}). More precisely, we want to find functions~$\tilde{\psi}^c, \tilde{\psi}^b : \mathbb{R}^{\Lambda_2^{(n)}} \to \mathbb{R}$ such that
\begin{align}
&\bar{H}_{\Lambda_2^{(n)}}\left( x^{\Lambda_2 ^{(n)}}\right) = \tilde{\psi}^c\left( x^{\Lambda_2 ^{(n)}} \right) + \tilde{\psi}^b \left( x^{\Lambda_2 ^{(n)}} \right), \label{e_psi_properties1} \\
&\Hess_{\mathbb{R}^{\Lambda_2^{(n)}}} \tilde{\psi}^c \geq c > 0 \qquad \text{ and } \qquad |\tilde{\psi}^b|_{\infty} \leq C < \infty. \label{e_psi_properties2}
\end{align}
Let us introduce auxiliary set~$E_n$ and Hamiltonian~$H_{\text{aux}} : \mathbb{R}^ N \to \mathbb{R}$. These are
\begin{align}
&E_n : = \{ k \in \Lambda  : \ \exists \ i  \in \Lambda_2 ^{(n)} \text{ such that } |k-i| \leq S    \}, \label{e_def_en}\\
&H_{\text{aux}}(x) := H(x) - \sum_{j \in E_n } \psi^b (x_j), \label{e_def_haux}
\end{align}
where~$S$ is a positive integer that will be chosen later. By definition, $H_{\text{aux}}$ is strictly convex in the space restricted to the spins~$x_i$ with~$i \in E_n$. The associated gce~$\mu_{\text{aux}}^{\sigma}(dx^{\Lambda_1})$ and the ce~$\mu_{\text{aux}, \widetilde{m}} \left(dx^{\Lambda_1}\right)$ are
\begin{align}
\mu^{\sigma}_{\text{aux}} \left(dx^{\Lambda_1}\right) = \mu^{\sigma}_{\text{aux}} \left(dx^{\Lambda_1} \big| x^{\Lambda_2} \right)  &: = \frac{1}{Z} \exp\left( \sigma \sum_{k \in \Lambda_1} x_k - H_{\text{aux}} \left(x^{\Lambda_1}, x^{\Lambda_2}\right) \right) dx^{\Lambda_1},
\end{align}
\begin{align}
\mu_{\text{aux}, \widetilde{m}} \left(dx^{\Lambda_1}\right) &= \mu_{\text{aux}, \widetilde{m}} \left(dx^{\Lambda_1} \big| x^{\Lambda_2}\right) \\ &:= \frac{1}{Z}\mathds{1}_{ \left\{ \frac{1}{K} \sum_{k \in \Lambda_1 } x_k = \widetilde{m} \right\}} \left(x^{\Lambda_1}\right)\exp\left( -H_{\text{aux}} \left(x^{\Lambda_1}, x^{\Lambda_2}\right)\right) \mathcal{L}^{K-1}(dx^{\Lambda_1}). \label{d_def_mu_aux,m}
\end{align}
As described in Section~\ref{s_auxiliary_lemmas}, we choose~$\sigma= \sigma( \widetilde{m})$ such that (cf. Definition~\ref{a_relation_sigma_m})
\begin{align} \label{e_m_and_mk}
\widetilde{m} = \frac{1}{K} \sum_{k \in \Lambda_1} \int  x_k  \mu^{\sigma}_{\text{aux}} (dx^{\Lambda_1}) =:  \frac{1}{K} \sum_{k \in \Lambda_1} \widetilde{m}_k . 
\end{align}
Motivated by the approach of~\cite{Me13}, the first attempt to decompose~$\bar{H}_{\Lambda_2^{(n)}}\left(x^{\Lambda_2^{(n)}} \right)$ according to~\eqref{e_psi_properties1} is to set
\begin{align}
\bar{H}_{\Lambda_2^{(n)}}\left(x^{\Lambda_2^{(n)}} \right) &= - \ln \int_{\frac{1}{K} \sum_{k \in \Lambda_1} x_k = \widetilde{m}} \exp \left( - H_{\text{aux}} \left(x^{\Lambda_1}, x^{\Lambda_2}\right) \right) \mathcal{L}(dx^{\Lambda_1}) \\
& \quad - \ln \frac{ \int_{\frac{1}{K} \sum_{k \in \Lambda_1} x_k = \widetilde{m}} \exp \left( - H \left(x^{\Lambda_1}, x^{\Lambda_2}\right) \right) \mathcal{L}(dx^{\Lambda_1}) }{\int_{\frac{1}{K} \sum_{k \in \Lambda_1} x_k = \widetilde{m}} \exp \left( - H_{\text{aux}} \left(x^{\Lambda_1}, x^{\Lambda_2}\right) \right) \mathcal{L}(dx^{\Lambda_1})}. \label{e_OR_naive_splitting}
\end{align}
The motivation for decomposition is the following. Because~$H_{\text{aux}}$ is strictly convex on the spins in~$E_n$, and log-integration conserves convexity, one expects that the first term on the right hand side of~\eqref{e_OR_naive_splitting} is strictly convex. One expects that the second term on the right hand side of~\eqref{e_OR_naive_splitting} is bounded because the Hamiltonian~$H$ and the auxiliary Hamiltonian~$H_{\text{aux}}$ are close up to bounded perturbations.   \\

However when doing so, it is not clear if the strict convexity of the first term on the right hand side of~\eqref{e_OR_naive_splitting} is uniform in the conditioned spins~$x^{\Lambda_2^{(l)}}$,~$l \neq n$. We circumvent this obstacle by adding and subtracting an additional term, i.e.~ 
\begin{align}
\bar{H}_{\Lambda_2^{(n)}}\left(x^{\Lambda_2^{(n)}} \right) & = - \left( \ln \int_{\frac{1}{K} \sum_{k \in \Lambda_1} x_k = \widetilde{m}} \exp \left( - H_{\text{aux}} \left(x^{\Lambda_1}, x^{\Lambda_2}\right) \right) \mathcal{L}(dx^{\Lambda_1})  \right.    \\
&\left. \qquad \quad + K \sigma(\widetilde{m})\widetilde{m}- \ln \int_{\mathbb{R}^{\Lambda_1}} \exp \left(\sigma(\widetilde{m})\sum_{k \in \Lambda_1} x_k - H_{\text{aux}} \left(x^{\Lambda_1}, x^{\Lambda_2}\right) \right) dx^{\Lambda_1}\right) \\
& \quad + \left( K \sigma(\widetilde{m})\widetilde{m}- \ln \int_{\mathbb{R}^{\Lambda_1}} \exp \left(\sigma(\widetilde{m})\sum_{k \in \Lambda_1} x_k - H_{\text{aux}} \left(x^{\Lambda_1}, x^{\Lambda_2}\right) \right) dx^{\Lambda_1} \right)  \\
& \quad - \ln \frac{ \int_{\frac{1}{K} \sum_{k \in \Lambda_1} x_k = \widetilde{m}} \exp \left( - H \left(x^{\Lambda_1}, x^{\Lambda_2}\right) \right) \mathcal{L}(dx^{\Lambda_1}) }{\int_{\frac{1}{K} \sum_{k \in \Lambda_1} x_k = \widetilde{m}} \exp \left( - H_{\text{aux}} \left(x^{\Lambda_1}, x^{\Lambda_2}\right) \right) \mathcal{L}(dx^{\Lambda_1})}. 
\end{align}
We define
\begin{align}
\tilde{\psi}^c \left(x^{\Lambda_2 ^{(n)}} \right) : &= K \sigma(\widetilde{m})\widetilde{m}- \ln \int_{\mathbb{R}^{\Lambda_1}} \exp \left(\sigma(\widetilde{m})\sum_{k \in \Lambda_1} x_k - H_{\text{aux}} \left(x^{\Lambda_1}, x^{\Lambda_2}\right) \right) dx^{\Lambda_1} , \label{e_aux_cramer_strict_convex} \\
\tilde{\psi}^b \left(x^{\Lambda_2 ^{(n)}} \right) : & = - \left( \ln \int_{\frac{1}{K} \sum_{k \in \Lambda_1} x_k = \widetilde{m}} \exp \left( - H_{\text{aux}} \left(x^{\Lambda_1}, x^{\Lambda_2}\right) \right) \mathcal{L}(dx^{\Lambda_1})  \right.    \\
&\left. + K \sigma(\widetilde{m})\widetilde{m}- \ln \int_{\mathbb{R}^{\Lambda_1}} \exp \left(\sigma(\widetilde{m})\sum_{k \in \Lambda_1} x_k - H_{\text{aux}} \left(x^{\Lambda_1}, x^{\Lambda_2}\right) \right) dx^{\Lambda_1}\right) \label{e_aux_cramer_representation}\\
& \quad - \ln \frac{ \int_{\frac{1}{K} \sum_{k \in \Lambda_1} x_k = \widetilde{m}} \exp \left( - H \left(x^{\Lambda_1}, x^{\Lambda_2}\right) \right) \mathcal{L}(dx^{\Lambda_1}) }{\int_{\frac{1}{K} \sum_{k \in \Lambda_1} x_k = \widetilde{m}} \exp \left( - H_{\text{aux}} \left(x^{\Lambda_1}, x^{\Lambda_2}\right) \right) \mathcal{L}(dx^{\Lambda_1})}. \label{e_aux_bounded_remaining_term}
\end{align}
Hence, we have the splitting
\begin{align}
\bar{H}_{\Lambda_2^{(n)}}\left(x^{\Lambda_2^{(n)}} \right)=    \tilde{\psi}^c \left(x^{\Lambda_2 ^{(n)}} \right)+ \tilde{\psi}^b \left(x^{\Lambda_2 ^{(n)}} \right).
\end{align}
Let us now explain, why the terms~$\tilde{\psi}^c$ and~$\tilde{\psi}^b$ heuristically satisfy~\eqref{e_psi_properties1} and~\eqref{e_psi_properties2}. Let us start with the property~\eqref{e_psi_properties2} i.e.~that~$\tilde{\psi}^b \left(x^{\Lambda_2 ^{(n)}} \right)$ is bounded. We observe that~$\tilde{\psi}^b \left(x^{\Lambda_2 ^{(n)}} \right)$ consists out of two terms, namely~\eqref{e_aux_cramer_representation} and~\eqref{e_aux_bounded_remaining_term}. The term~\eqref{e_aux_bounded_remaining_term} is unchanged from the first attempt and should be bounded by the same argument as above. The term~\eqref{e_aux_cramer_representation} is representing the difference between the free energy of an canonical ensemble and the free energy of an associated modified grand-canonical ensemble. Hence, this term should be uniformly bounded by the equivalence of ensembles.\newline
Let us explain why~\eqref{e_psi_properties1} is satisfied i.e.~why the term~$\tilde{\psi}^c$ is uniformly strictly convex. We observe that~$\tilde{\psi}^c$ is a Legendre-transformation. By basic properties of the Legendre transform the convexity properties of $\tilde{\psi}^c$ are determined by the convexity properties of the log-partition function of the associated grand-canonical ensemble. In the first attempt we would have to analyze the convexity of the log-partition function of a canonical ensemble. Now, we only have to analyze the convexity of the log-partition function of a grand-canonical ensemble. This is easier to analyze and turns out to be more robust wrt.~the conditioned spin-values.\\

This argument is made rigorous in the following two statements.

\begin{lemma} \label{l_strictly_convex}
For both~$K=|\Lambda_1|$ and~$S$ large enough (cf.~\eqref{e_def_en}), the function~$\tilde{\psi}^c$ is strictly convex in the sense that there is a positive constant~$c>0$ with
\begin{align}
\Hess_{\mathbb{R}^{\Lambda_2 ^{(n)}}} \tilde{\psi}^c \geq c \Id_{\mathbb{R}^{\Lambda_2 ^{(n)}}},
\end{align}
where~$\Id_{\mathbb{R}^{\Lambda_2 ^{(n)}}}$ denotes the identity map on~$\mathbb{R}^{\Lambda_2 ^{(n)}}$. 
\end{lemma}

\begin{lemma} \label{l_uniformly_bounded}
The function~$\tilde{\psi}^b$ is uniformly bounded. More precisely, it holds that
\begin{align}
| \tilde{\psi}^b |_{\infty} \lesssim 1 + R + S.
\end{align}
\end{lemma}
We shall see how the statements from above yield Lemma~\ref{l_block_marginal_lsi}. \\

\noindent \emph{Proof of Lemma~\ref{l_block_marginal_lsi}.} \ 
By choosing~$S$ large but fixed, the block Hamiltonian~$\bar{H}_{\Lambda_2^{(n)}}$ is a sum of strictly convex function~$\widetilde{\psi}^c$ and a bounded function~$\widetilde{\psi}^b$. Then the Lemma~\ref{l_block_marginal_lsi} follows from a combination of Bakry-\'{E}mery criterion (see Theorem~\ref{a_bakry_emery}) and Holley-Stroock Perturbation Principle (see Theorem~\ref{a_holley_stroock}). \qed

\medskip

It remains to prove Lemma~\ref{l_strictly_convex} and Lemma~\ref{l_uniformly_bounded}. Let us begin with providing auxiliary statements which will be verified in Section~\ref{s_strict_convexity_aux_proof}. For the notational convenience, we simply write~$\tilde{\sigma}= \sigma(\widetilde{m})$ (cf. see~\eqref{e_m_and_mk}). \\

For the proof of Lemma~\ref{l_strictly_convex} we need to show that there exists a positive constant~$c>0$ such that the following holds for~$K$ and~$S$ large enough.
\begin{align}
\Hess_{\mathbb{R}^{\Lambda_2 ^{(n)}}} \left( K \tilde{\sigma} \widetilde{m}- \ln \int_{\mathbb{R}^{ \Lambda_1 }} \exp \left(\tilde{\sigma} \sum_{k \in \Lambda_1} x_k - H_{\text{aux}} \left(x^{\Lambda_1}, x^{\Lambda_2}\right) \right) dx^{\Lambda_1} \right)\ \geq c \Id_{\mathbb{R}^{ \Lambda_2^{(n)} }}, \label{e_claim_hess_pos}
\end{align}
where~$\Id_{\mathbb{R}^{\Lambda_2 ^{(n)}}}$ denotes the identity map on~$\mathbb{R}^{\Lambda_2^{(n)}}$. As a first step, we calculate a formula for the left hand side of~\eqref{e_claim_hess_pos}.

\begin{lemma}\label{l_second_expression} For each~$i, j \in \Lambda_2^{(n)}$, it holds that
\begin{align}
&\frac{d^2}{dx_i dx_j } \left( K \tilde{\sigma}\widetilde{m}- \ln \int_{\mathbb{R}^{\Lambda_1}} \exp \left(\tilde{\sigma} \sum_{k \in \Lambda_1} x_k - H_{\text{aux}} \left(x^{\Lambda_1}, x^{\Lambda_2}\right) \right) dx^{\Lambda_1} \right) \\
& \qquad = \var_{\mu_{\text{aux}}^{\tilde{\sigma}}} \left( \sum_{k \in \Lambda_1} X_k \right) \frac{d \tilde{\sigma}}{dx_i} \frac{d \tilde{\sigma}}{dx_j} \label{e_hessian_variance_term} \\
& \qquad \quad + \mathbb{E}_{\mu_{\text{aux}}^{\tilde{\sigma}}} \left[   \frac{\partial ^2}{\partial x_i \partial x_j} H_{\text{aux}} (X)\right] - \cov_{\mu_{\text{aux}}^{\tilde{\sigma}}} \left(\frac{\partial}{\partial x_i} H_{\text{aux}}(X), \frac{\partial}{\partial x_j} H_{\text{aux}}(X) \right). \label{e_hessian_covariance_term}
\end{align}
\end{lemma}

\medskip

Next, we provide an auxiliary estimate of partial derivative of~$\tilde{\sigma}$ with respect to~$x_i$, where~$i \in \Lambda_2^{(n)}$.

\begin{lemma}\label{l_dsigma_dxi_bounds}
For each~$i \in \Lambda_2 ^{(n)}$, it holds that
\begin{align}
\left| \frac{d \tilde{\sigma}}{dx_i} \right| \lesssim \frac{1}{K}. \label{e_first_der_sigma_xi}
\end{align}
\end{lemma}

\medskip

The last step towards to the proof of Lemma~\ref{l_strictly_convex} is the strict convexity of~\eqref{e_hessian_covariance_term}.

\begin{lemma} \label{l_Brascamp_lieb}
There is a positive constant~$c>0$ such that for both~$K=|\Lambda_1|$ and~$S$ large enough (cf.~\eqref{e_def_en}) the following holds.
\begin{align}
\left(\mathbb{E}_{\mu_{\text{aux}}^{\tilde{\sigma}}} \left[ \frac{\partial ^2}{\partial x_i \partial x_j } H_{\text{aux}}(X) \right] - \cov_{\mu_{\text{aux}}^{\tilde{\sigma}}} \left( \frac{\partial}{\partial x_i }H_{\text{aux}}(X), \frac{\partial}{\partial x_j }H_{\text{aux}}(X) \right)\right)_{i, j \in \Lambda_2 ^{(n)}} \geq c \Id_{\mathbb{R}^{ \Lambda_2 ^{(n)} }}. \label{e_lemma_brascamp_lieb_statement}
\end{align}
\end{lemma}

\medskip

The proof of Lemma~\ref{l_second_expression}, Lemma~\ref{l_dsigma_dxi_bounds} and Lemma~\ref{l_Brascamp_lieb} are presented in Section~\ref{s_strict_convexity_aux_proof}. Let us see how these statements yield Lemma~\ref{l_strictly_convex}. \\

\noindent  \emph{Proof of Lemma~\ref{l_strictly_convex}.} \ Due to Lemma~\ref{l_dsigma_dxi_bounds} it holds that for~$i, j \in \Lambda_2^{(n)}$
\begin{align}
&\frac{d^2}{dx_i dx_j } \left( K \tilde{\sigma} \widetilde{m}- \ln \int_{\mathbb{R}^{\Lambda_1}} \exp \left(\tilde{\sigma}\sum_{k \in \Lambda_1} x_k - H_{\text{aux}} \left(x^{\Lambda_1}, x^{\Lambda_2}\right) \right) dx^{\Lambda_1} \right) = T_{\eqref{e_hessian_variance_term}} + T_{\eqref{e_hessian_covariance_term}}. 
\end{align}
Hence, the desired statement follows from Lemma~\ref{l_Brascamp_lieb} and the estimate
\begin{align}
 \left| T_{\eqref{e_hessian_variance_term}} \right| \lesssim \frac{1}{K}. \label{e_bound_on_hessian_variance_term}
\end{align}
Indeed,~\eqref{e_bound_on_hessian_variance_term} directly follows from a combination of Lemma~\ref{l_dsigma_dxi_bounds} and the observation that (cf. Lemma~\ref{l_variance_estimate})
\begin{align}
\var_{\mu_{\text{aux}}^{\tilde{\sigma}}} \left( \sum_{k \in \Lambda_1} X_k \right) \lesssim K.
\end{align}
\qed 

\medskip

Let us turn to the proof of Lemma~\ref{l_uniformly_bounded}. We need an auxiliary statement to begin with. \\
\begin{lemma} \label{l_aux_cramers_representation}
Let~$X= (X_k)_{k \in \Lambda_1}$ be a random variable distributed according to the measure~$\mu_{\text{aux}}^{\tilde{\sigma}}$. Denote $g_{\text{aux}}$ by the density of the random variable~$\frac{1}{\sqrt{ K}} \sum_{k \in \Lambda_1} \left( X_k - \widetilde{m}_k \right)$. It holds that
\begin{align}
\ln g_{\text{aux}}(0) &= \left( \ln \int_{\frac{1}{K} \sum_{k \in \Lambda_1} x_k = \widetilde{m}} \exp \left( - H_{\text{aux}} \left(x^{\Lambda_1}, x^{\Lambda_2}\right) \right) \mathcal{L}(dx^{\Lambda_1})  \right.    \\
&\left. \qquad \quad + K \tilde{\sigma}\widetilde{m}- \ln \int_{\mathbb{R}^{\Lambda_1}} \exp \left(\tilde{\sigma}\sum_{k \in \Lambda_1} x_k - H_{\text{aux}} \left(x^{\Lambda_1}, x^{\Lambda_2}\right) \right) dx^{\Lambda_1}\right).
\end{align}
\end{lemma}

\medskip The proof of Lemma~\ref{l_aux_cramers_representation} is given in Section~\ref{s_strict_convexity_aux_proof}. Let us now turn to the proof of Lemma~\ref{l_uniformly_bounded}. \\

\noindent \emph{Proof of Lemma~\ref{l_uniformly_bounded}.} \ 
We observe that 
\begin{align}
\tilde{\psi}^b \left(x^{\Lambda_2 ^{(n)}} \right)  = T_{\eqref{e_aux_cramer_representation}} + T_{\eqref{e_aux_bounded_remaining_term}}.
\end{align}
Estimation of~$T_{\eqref{e_aux_cramer_representation}}$: Because~$\mu_{\text{aux}}^{\sigma}$ is a one-dimensional gce, the desired estimate
\begin{align} \label{e_estimation_cramer_rep}
|T_{\eqref{e_aux_cramer_representation}}| \lesssim 1
\end{align}
follows from a combination of Lemma~\ref{l_aux_cramers_representation} and Proposition~\ref{p_main_computation}. \\

\noindent Estimation of~$T_{\eqref{e_aux_bounded_remaining_term}}$: It holds that
\begin{align}
\left|T_{\eqref{e_aux_bounded_remaining_term}} \right| &= \left|\ln \frac{ \int_{\frac{1}{K} \sum_{k \in \Lambda_1} x_k = \widetilde{m}} \exp \left( - H \left(x^{\Lambda_1}, x^{\Lambda_2}\right) \right) \mathcal{L}(dx^{\Lambda_1}) }{\int_{\frac{1}{K} \sum_{k \in \Lambda_1} x_k = \widetilde{m}} \exp \left( - H_{\text{aux}} \left(x^{\Lambda_1}, x^{\Lambda_2}\right) \right) \mathcal{L}(dx^{\Lambda_1})} \right| \\
& \overset{\eqref{d_def_mu_aux,m}}{=}\left| \ln \int \exp\left( -\sum_{ j \in E_n} \psi^{b} (x_j ) \right)\mu_{\text{aux}, \widetilde{m}}\left(dx^{\Lambda_1}\right)   \right|  \\
& \leq \sum_{j \in E_n } \| \psi^{b} \|_{\infty} \lesssim R+S. \label{e_estimation_aux_bounded_remaining}
\end{align}
Then a combination of~\eqref{e_estimation_cramer_rep} and~\eqref{e_estimation_aux_bounded_remaining} yields, as desired,
\begin{align}
|\tilde{\psi}^b|_{\infty} \lesssim 1 + R + S.
\end{align}
\qed

\medskip

\subsection{Proof of Auxiliary Statements} \label{s_strict_convexity_aux_proof}

In this section, we provide the proof of Lemma~\ref{l_second_expression}, Lemma~\ref{l_dsigma_dxi_bounds}, Lemma~\ref{l_Brascamp_lieb} and Lemma~\ref{l_aux_cramers_representation}. Let us begin with providing the proof of Lemma~\ref{l_second_expression}. \\

\noindent  \emph{Proof of Lemma~\ref{l_second_expression}.} \ 
A direct calculation yields that
\begin{align}
&\frac{d}{dx_j } \left( K \tilde{\sigma}\widetilde{m}- \ln \int_{\mathbb{R}^{\Lambda_1}} \exp \left(\tilde{\sigma}\sum_{k \in \Lambda_1} x_k - H_{\text{aux}} \left(x^{\Lambda_1}, x^{\Lambda_2}\right) \right) dx^{\Lambda_1} \right) \\
&\qquad = K \left(\frac{d}{dx_j}\tilde{\sigma}\right) \widetilde{m} + K \tilde{\sigma} \frac{d \widetilde{m}}{dx_j} \\
&\qquad \quad - \frac{\partial}{\partial \tilde{\sigma}}\left( \ln \int_{\mathbb{R}^{\Lambda_1}} \exp \left(\tilde{\sigma}\sum_{k \in \Lambda_1} x_k - H_{\text{aux}} \left(x^{\Lambda_1}, x^{\Lambda_2}\right) \right) dx^{\Lambda_1} \right) \frac{d\tilde{\sigma}}{dx_j} \\
&\qquad \quad - \frac{\partial}{\partial x_j}\left( \ln \int_{\mathbb{R}^{\Lambda_1}} \exp \left(\tilde{\sigma}\sum_{k \in \Lambda_1} x_k - H_{\text{aux}} \left(x^{\Lambda_1}, x^{\Lambda_2}\right) \right) dx^{\Lambda_1} \right) \\
&\qquad \overset{\eqref{e_m_and_mk}}{=}  K \tilde{\sigma} \frac{d \widetilde{m}}{dx_j} - \frac{\partial}{\partial x_j}\left( \ln \int_{\mathbb{R}^{\Lambda_1}} \exp \left(\tilde{\sigma}\sum_{k \in \Lambda_1} x_k - H_{\text{aux}} \left(x^{\Lambda_1}, x^{\Lambda_2}\right) \right) dx^{\Lambda_1} \right) \\
&\qquad \overset{\eqref{e_def_tilde_m}}{=} - \tilde{\sigma} - \mathbb{E}_{\mu_{\text{aux}}^{\tilde{\sigma}}} \left[  - \frac{\partial}{\partial x_j} H_{\text{aux}} (X)\right] \\
&\qquad = - \tilde{\sigma} + \mathbb{E}_{\mu_{\text{aux}}^{\tilde{\sigma}}} \left[   \frac{\partial}{\partial x_j} H_{\text{aux}} (X)\right].
\end{align}
Taking further derivative, one gets
\begin{align}
&\frac{d^2}{dx_i dx_j } \left( K \tilde{\sigma}\widetilde{m}- \ln \int_{\mathbb{R}^{\Lambda_1}} \exp \left(\tilde{\sigma}\sum_{k \in \Lambda_1} x_k - H_{\text{aux}} \left(x^{\Lambda_1}, x^{\Lambda_2}\right) \right) dx^{\Lambda_1} \right) \\
& \qquad = - \frac{d \tilde{\sigma}}{dx_i} + \frac{d}{dx_i} \mathbb{E}_{\mu_{\text{aux}}^{\tilde{\sigma}}} \left[   \frac{\partial}{\partial x_j} H_{\text{aux}} (X)\right] \\
& \qquad = -  \frac{d \tilde{\sigma}}{dx_i} +  \frac{\partial}{\partial \tilde{\sigma}} \mathbb{E}_{\mu_{\text{aux}}^{\tilde{\sigma}}} \left[   \frac{\partial}{\partial x_j} H_{\text{aux}} (X)\right]\frac{d \tilde{\sigma}}{dx_i} +  \frac{\partial}{\partial x_i} \mathbb{E}_{\mu_{\text{aux}}^{\tilde{\sigma}}} \left[   \frac{\partial}{\partial x_j} H_{\text{aux}} (X)\right]   \\
& \qquad =  \frac{d \tilde{\sigma}}{dx_i} \left( \cov_{\mu_{\text{aux}}^{\tilde{\sigma}}} \left( \sum_{k \in \Lambda_1}  X_k  , \frac{\partial}{\partial x_j}H_{\text{aux}}(X)  \right) -1 \right) \\
& \qquad \quad + \mathbb{E}_{\mu_{\text{aux}}^{\tilde{\sigma}}} \left[   \frac{\partial^2}{\partial x_i \partial x_j} H_{\text{aux}} (X)\right] - \cov_{\mu_{\text{aux}}^{\tilde{\sigma}}} \left(\frac{\partial}{\partial x_i} H_{\text{aux}}(X), \frac{\partial}{\partial x_j} H_{\text{aux}}(X) \right) \\
& \qquad \overset{\eqref{e_sigma_first_derivative_wrt_xi}}{=} \var_{\mu_{\text{aux}}^{\tilde{\sigma}}} \left( \sum_{k \in \Lambda_1} X_k \right) \frac{d \tilde{\sigma}}{dx_i} \frac{d \tilde{\sigma}}{dx_j} \\
& \qquad \quad + \mathbb{E}_{\mu_{\text{aux}}^{\tilde{\sigma}}} \left[   \frac{\partial ^2}{\partial x_i \partial x_j} H_{\text{aux}} (X)\right] - \cov_{\mu_{\text{aux}}^{\tilde{\sigma}}} \left(\frac{\partial}{\partial x_i} H_{\text{aux}}(X), \frac{\partial}{\partial x_j} H_{\text{aux}}(X) \right).  \label{e_dxidxj_remaining}
\end{align}
\qed

Let us now turn to the proof of Lemma~\ref{l_dsigma_dxi_bounds}. \\

\noindent \emph{Proof of Lemma~\ref{l_dsigma_dxi_bounds}.} \ A combination of the definition~\eqref{e_def_tilde_m} of~$\widetilde{m}$ and the equality~\eqref{e_m_and_mk} yields
\begin{align}
\frac{N}{K} m - \frac{1}{K} \sum_{i \in \Lambda_2} x_i = \frac{1}{K} \sum_{k \in \Lambda_1} \widetilde{m}_k.
\end{align}
We recall that~$\mu_{\text{aux}}^{\tilde{\sigma}} \left(dx^{\Lambda_1}\right) = \mu_{\text{aux}}^{\tilde{\sigma}} \left(dx^{\Lambda_1} \big| x^{\Lambda_2}\right)$ depends on~$x_i, i \in \Lambda_2^{(n)}$ via the field~$\tilde{\sigma}$ and the conditioned spin value~$x^{\Lambda_2}$. Then differentiating both sides with respect to~$x_i, i \in \Lambda_2^{(n)}$ gives
\begin{align}
- \frac{1}{K} &= \frac{1}{K}  \frac{d}{dx_i}\left( \sum_{k \in \Lambda_1} \widetilde{m}_k \right) \\
& = \frac{1}{K} \frac{\partial}{\partial \tilde{\sigma}} \left(\sum_{k \in \Lambda_1}  \widetilde{m}_k\right) \frac{d\tilde{\sigma}}{dx_i} + \frac{1}{K}\frac{\partial}{\partial x_i} \left( \sum_{k \in \Lambda_1} \widetilde{m}_k  \right) \\
& =  \frac{1}{K}  \var_{\mu_{\text{aux}}^{\tilde{\sigma}}} \left( \sum_{k \in \Lambda_1} X_k \right) \frac{d\tilde{\sigma}}{dx_i} + \frac{1}{K} \cov_{\mu_{\text{aux}}^{\tilde{\sigma}}} \left( \sum_{k \in \Lambda_1} X_k, -\frac{\partial}{\partial x_i} H_{\text{aux}} (X) \right). \label{e_sigma_first_derivative_wrt_xi}
\end{align}
Our strategy is to estimate every factor appearing in~\eqref{e_sigma_first_derivative_wrt_xi}. Then rearranging~\eqref{e_sigma_first_derivative_wrt_xi} will yield the desired estimate. Note that the measure~$\mu_{\text{aux}}^{\tilde{\sigma}} $ is only active on~$\Lambda_1$ while~$i \in \Lambda_2 ^{(n)}$. Since the covariance is invariant under adding constants, it holds that
\begin{align}
\left| \cov_{\mu_{\text{aux}}^{\tilde{\sigma}}} \left( \sum_{k \in \Lambda_1} X_k, \frac{\partial}{\partial x_i} H_{\text{aux}} (X) \right) \right|
&= \left| \cov_{\mu_{\text{aux}}^{\tilde{\sigma}}} \left( \sum_{k \in \Lambda_1} X_k, X_i  + s_i + \frac{1}{2} \sum_{ l : 1\leq |l-i| \leq R} M_{il} X_l \right) \right| \\
&=\frac{1}{2} \left| \cov_{\mu_{\text{aux}}^{\tilde{\sigma}}} \left( \sum_{k \in \Lambda_1} X_k,  \sum_{ \substack{l \in \Lambda_1 \\ l : 1\leq |l-i| \leq R}} M_{il} X_l \right) \right|. \label{e_aux_exp_app}
\end{align}
We note that the properties described in Section~\ref{s_auxiliary_lemmas} hold for~$\mu_{\text{aux}}^{\tilde{\sigma}}$ because it is a gce in the one-dimensional lattice. In particular~$\mu_{\text{aux}}^{\tilde{\sigma}}$ has exponential decay of correlations as in Lemma~\ref{l_exponential_decay_of_correlations}. Therefore it holds that
\begin{align}
T_{\eqref{e_aux_exp_app}} \leq \frac{1}{2} \sum_{k \in \Lambda_1} \sum_{ \substack{l \in \Lambda_1 \\ l : 1\leq |l-i| \leq R}} \left| \cov_{\mu_{\text{aux}}^{\tilde{\sigma}}} \left(  X_k,   M_{il} X_l \right) \right| \lesssim  R \lesssim 1. \label{e_estimate_cov}
\end{align}
As a consequence, by rearranging~\eqref{e_sigma_first_derivative_wrt_xi} we get
\begin{align}
\left|\frac{d\tilde{\sigma}}{dx_i} \right| = \frac{K}{\var_{\mu_{\text{aux}}^{\tilde{\sigma}}}\left( \sum_{k \in \Lambda_1 } X_k \right)  }\left| \frac{1}{K} \cov_{\mu_{\text{aux}}^{\tilde{\sigma}}} \left( \sum_{k \in \Lambda_1} X_k, \frac{\partial}{\partial x_i} H_{\text{aux}} (X) \right) - \frac{1}{K}  \right| \overset{Lemma~\ref{l_variance_estimate},~\eqref{e_estimate_cov}}{\lesssim} \frac{1}{K}.
\end{align}
\qed

The proof of Lemma~\ref{l_Brascamp_lieb} consists of a lengthy calculations. We provide the proof of Lemma~\ref{l_Brascamp_lieb} modulo an auxiliary statement (see Lemma~\ref{l_aux_lem_for_brascamp_lieb} from below). \\

\noindent \emph{Proof of Lemma~\ref{l_Brascamp_lieb}.} \ The argument is inspired by calculations done in~\cite[Lemma 3.1]{Me13}. The main difference is that we consider the multi-variable case (in $\mathbb{R}^{\Lambda_2^{(n)}}$) while~\cite[Lemma 3.1]{Me13} only considers the single variable case. Recall the definition~\eqref{e_def_en} of the set~$E_n$. Let us decompose the auxiliary measure~$\mu_{\text{aux}}^{\tilde{\sigma}}$ into the conditional and the marginal measure as follows:
\begin{align}
\mu_{\text{aux}}^{\tilde{\sigma}} \left( dx^{\Lambda_1} \right) = \mu_{\text{aux}}^{\tilde{\sigma}} \left( \left( dx_k \right)_{ k \in \Lambda_1,  k \in E_n} \big| \left( x_l \right)_{l \in \Lambda_1, l \not\in E_n } \right) \bar{\mu}_{\text{aux}}^{\tilde{\sigma}} \left( \left(dx_l\right)_{l \in \Lambda_1 , l \not\in E_n}  \right).
\end{align}
We write~$\mu_{\text{aux,c}}^{\tilde{\sigma}} = \mu_{\text{aux}}^{\tilde{\sigma}} \left( \left( dx_k \right)_{ k \in \Lambda_1, k \in E_n} \big| \left( x_l \right)_{l \in \Lambda_1, l \not\in E_n} \right)$ for notational convenience. Then it follows that 

\begin{align}
&\mathbb{E}_{\mu_{\text{aux}}^{\tilde{\sigma}}} \left[ \frac{\partial^2}{\partial x_i \partial x_j }H_{\text{aux}}(X) \right] - \cov_{\mu_{\text{aux}}^{\tilde{\sigma}}}\left( \frac{\partial}{\partial x_i }H_{\text{aux}}(X), \frac{\partial}{\partial x_j }H_{\text{aux}}(X)\right) \\
& \qquad = \int \left( \int \frac{\partial^2}{\partial x_i \partial x_j}H_{\text{aux}}(x) \mu_{\text{aux,c}}^{\tilde{\sigma}} - \cov_{\mu_{\text{aux,c}}} \left( \frac{\partial}{\partial x_i }H_{\text{aux}}(x), \frac{\partial}{\partial x_j }H_{\text{aux}}(x) \right)   \right)\bar{\mu}_{\text{aux}}^{\tilde{\sigma}} \\
& \qquad \quad - \cov_{\bar{\mu}_{\text{aux}}^{\tilde{\sigma}}}\left( \int \frac{\partial}{\partial x_i} H_{\text{aux}}(x)  \mu_{\text{aux,c}}^{\tilde{\sigma}} , \int \frac{\partial}{\partial x_j} H_{\text{aux}}(x) \mu_{\text{aux,c}}^{\tilde{\sigma}}  \right). \label{e_brascamp_decomposition}
\end{align}
As noted before, if restricted to spins~$x_i$ with~$i \in E_n$, the Hamiltonian~$H_{\text{aux}}$ is strictly convex. Then by Brascamp-Lieb inequality (cf.~\cite{Bra76}), there is a constant~$c>0$ with
\begin{align}
\left( \int \frac{\partial^2}{\partial x_i \partial x_j }H_{\text{aux}}(x) \mu_{\text{aux,c}}^{\tilde{\sigma}} - \cov_{\mu_{\text{aux,c}}^{\tilde{\sigma}}} \left( \frac{\partial}{\partial x_i }H_{\text{aux}}(x), \frac{\partial}{\partial x_j }H_{\text{aux}}(x) \right) \right)_{i,j} \geq c\Id_{\mathbb{R}^{ \Lambda_2 ^{(n)}}}.
\end{align}
and as a consequence,
\begin{align}
\left(\int \left( \int \frac{\partial^2}{\partial x_i \partial x_j}H_{\text{aux}}(x) \mu_{\text{aux,c}}^{\tilde{\sigma}} - \cov_{\mu_{\text{aux,c}}^{\tilde{\sigma}}} \left( \frac{\partial}{\partial x_i }H_{\text{aux}} (x) , \frac{\partial}{\partial x_i }H_{\text{aux}} (x)\right)   \right) \bar{\mu}_{\text{aux}}^{\tilde{\sigma}} \right)_{i, j} \geq c\Id_{\mathbb{R}^{ \Lambda_2^{(n)} }}. \label{e_brascamp_lieb}
\end{align}
Then a combination of~\eqref{e_brascamp_decomposition},~\eqref{e_brascamp_lieb} and Lemma~\ref{l_aux_lem_for_brascamp_lieb} from below yields the desired property~\eqref{e_lemma_brascamp_lieb_statement} by choosing~$S$ large enough. 
\qed 

\medskip

The following statement, Lemma~\ref{l_aux_lem_for_brascamp_lieb}, is an estimation of the second therm in the right hand side of~\eqref{e_brascamp_decomposition}. \\

\begin{lemma} \label{l_aux_lem_for_brascamp_lieb}
Recall the definition~\eqref{e_def_en} of the set~$E_n$.
\begin{align}
E_n : = \{ k \in \Lambda  : \ \exists \ i  \in \Lambda_2 ^{(n)} \text{ such that } |k-i| \leq S    \}.
\end{align}
For each~$i, j \in \Lambda_2 ^{(n)}$, it holds that
\begin{align}
\left| \cov_{\bar{\mu}_{\text{aux}}^{\tilde{\sigma}}}\left( \int \frac{\partial}{\partial x_i} H_{\text{aux}}(x)  \mu_{\text{aux,c}}^{\tilde{\sigma}} , \int \frac{\partial}{\partial x_j} H_{\text{aux}}(x)  \mu_{\text{aux,c}}^{\tilde{\sigma}}  \right) \right| \lesssim S^2 \exp\left(-CS\right).
\end{align}
\end{lemma}

\noindent \emph{Proof of Lemma~\ref{l_aux_lem_for_brascamp_lieb}.} \ We begin with a simple observation:
\begin{align}
\frac{\partial}{\partial x_i} H_{\text{aux}}\left(x\right) = x_i + s_i + \frac{1}{2} \sum_{p : 1 \leq |p-i| \leq R} M_{ip} x_p .
\end{align}
As the measures~$\mu_{\text{aux,c}}^{\tilde{\sigma}}$ and~$\bar{\mu}_{\text{aux}}^{\tilde{\sigma}}$ are defined in the subspace of~$x^{\Lambda_1}$, we may regard~$x_l$'s with~$l \in \Lambda_2$ as constants. As a consequence, it holds that
\begin{align}
& \cov_{\bar{\mu}_{\text{aux}}^{\tilde{\sigma}}}\left( \int \frac{\partial}{\partial x_i} H_{\text{aux}}(x)  \mu_{\text{aux,c}}^{\tilde{\sigma}} , \int \frac{\partial}{\partial x_j} H_{\text{aux}}(x)  \mu_{\text{aux,c}}^{\tilde{\sigma}}  \right) \\
& \quad = \frac{1}{4} \cov_{\bar{\mu}_{\text{aux}}^{\tilde{\sigma}}}\left( \int \sum_{\substack{ p \in \Lambda_1, \\ 1\leq |p-i| \leq R }}  M_{ip} x_p \ \mu_{\text{aux,c}}^{\tilde{\sigma}}, \int \sum_{\substack{ q \in \Lambda_1, \\ 1\leq |q-j| \leq R }}  M_{jq} x_q \ \mu_{\text{aux,c}}^{\tilde{\sigma}}  \right).
\end{align}
To estimate the covariance from above, let us double the variables to get
\begin{align}
&\cov_{\bar{\mu}_{\text{aux}}^{\tilde{\sigma}}}\left(\int \sum_{\substack{ p \in \Lambda_1, \\ 1\leq |p-i| \leq R }}  M_{ip} x_p \ \mu_{\text{aux,c}}^{\tilde{\sigma}} , \int \sum_{\substack{ q \in \Lambda_1, \\ 1\leq |q-j| \leq R }}  M_{jq} x_q \ \mu_{\text{aux,c}}^{\tilde{\sigma}} \right) \\
&  = \int \int  \left( \int \sum_{\substack{ p \in \Lambda_1, \\ 1\leq |p-i| \leq R }}  M_{ip} x_p \ \mu_{\text{aux,c}}^{\tilde{\sigma}}(dx | y) - \int \sum_{\substack{ p \in \Lambda_1, \\ 1\leq |p-i| \leq R }}  M_{ip} x_p \ \mu_{\text{aux,c}}^{\tilde{\sigma}}(dx | z)  \right) \\
& \times \left( \int \sum_{\substack{ q \in \Lambda_1, \\ 1\leq |q-j| \leq R }}  M_{jq} x_q \ \mu_{\text{aux,c}}^{\tilde{\sigma}}(dx | y) - \int \sum_{\substack{ q \in \Lambda_1, \\ 1\leq |q-j| \leq R }}  M_{jq} x_q \ \mu_{\text{aux,c}}^{\tilde{\sigma}}(dx | z)  \right) \bar{\mu}_{\text{aux}}^{\tilde{\sigma}}(dy)\bar{\mu}_{\text{aux}}^{\tilde{\sigma}}(dz). \label{e_double_variable_covariance}
\end{align}
Then it follows from the fundamental theorem of calculus that
\begin{align}
&\left| \int \sum_{\substack{ p \in \Lambda_1, \\ 1\leq |p-i| \leq R }}  M_{ip} x_p \ \mu_{\text{aux,c}}^{\tilde{\sigma}}(dx | y) - \int \sum_{\substack{ p \in \Lambda_1, \\ 1\leq |p-i| \leq R }}  M_{ip} x_p \ \mu_{\text{aux,c}}^{\tilde{\sigma}}(dx | z) \right| \\
& \quad= \left| \int_0 ^1 \left( \frac{d}{dt} \int \sum_{\substack{ p \in \Lambda_1, \\ 1\leq |p-i| \leq R }}  M_{ip} x_p \ \mu_{\text{aux,c}}^{\tilde{\sigma}}(dx | ty + (1-t)z) \right)dt \right| \label{e_fundamental_thm_calc} \\
& \quad = \left|\int_0 ^1 \cov_{\mu_{\text{aux,c}}^{\tilde{\sigma}}(dx \mid ty+(1-t)z)} \left(  \sum_{\substack{ p \in \Lambda_1, \\ 1\leq|p-i| \leq R }}  M_{ip} x_p,  \sum_{\substack{ r_1 \in E_n \\ s_1 \notin E_n \\ 1 \leq |r_1-s_1| \leq R }} M_{r_1 s_1 }x_{r_1} ( y_{s_1} -z_{s_1} ) \right)  dt \right| \\
& \quad \leq \sum_{\substack{ p \in \Lambda_1, \\ 1\leq |p-i| \leq R }} \sum_{\substack{ r_1 \in E_n \\ s_1 \notin E_n \\ 1\leq |r_1-s_1| \leq R }}  |M_{ip}| |M_{r_1 s_1}| |y_{s_1} -z_{s_1}| \sup_{0 \leq t \leq 1} \left| \cov_{\mu_{\text{aux,c}}^{\tilde{\sigma}}(dx \mid ty+(1-t)z)} \left(x_p, x_{r_1}\right)\right|. \label{e_cov_aux_c_exp_decay}
\end{align}
We note that~$\mu_{\text{aux}}^{\tilde{\sigma}}$,~$\mu_{\text{aux,c}}^{\tilde{\sigma}}$ are also gces on the one-dimensional lattice and they satisfy properties listed in Section~\ref{s_auxiliary_lemmas}. Therefore an application of Lemma~\ref{l_exponential_decay_of_correlations} yields
\begin{align}
T_{\eqref{e_cov_aux_c_exp_decay}} \lesssim \sum_{\substack{ p \in \Lambda_1, \\ 1 \leq |p-i| \leq R }} \sum_{\substack{ r_1 \in E_n \\ s_1 \notin E_n \\ 1\leq |r_1-s_1| \leq R }} |M_{ip}| |M_{r_1 s_1}| |y_{s_1} -z_{s_1}| \exp\left(-C | p-r_1 | \right). \label{e_d_dt_int_aux1}
\end{align}
Similarly, one gets
\begin{align}
&\left| \int \sum_{\substack{ q \in \Lambda_1, \\ |q-j| \leq R }}  M_{ip} x_p \ \mu_{\text{aux,c}}^{\tilde{\sigma}}(dx | y) - \int \sum_{\substack{ q \in \Lambda_1, \\ |q-j| \leq R }}  M_{jq} x_q \ \mu_{\text{aux,c}}^{\tilde{\sigma}}(dx | z) \right| \\
& \qquad \lesssim \sum_{\substack{ q \in \Lambda_1, \\ |q-j| \leq R }} \sum_{\substack{ r_2 \in E_n \\ s_2 \notin E_n \\ |r_2-s_2| \leq R }} |M_{jq}| |M_{r_2 s_2}| |y_{s_2} -z_{s_2}| \exp\left(-C | q-r_2 | \right). \label{e_d_dt_int_aux2}
\end{align}
Note that
\begin{align}
\int \int \left| y_s -z_s \right| \left| y_r -z_r \right| \bar{\mu}_{\text{aux}}^{\tilde{\sigma}} (dy)\bar{\mu}_{\text{aux}}^{\tilde{\sigma}} (dz) & \leq \frac{1}{2} \int \int \left( y_s -z_s \right)^2 + \left( y_r -z_r \right)^2 \bar{\mu}_{\text{aux}}^{\tilde{\sigma}} (dy)\bar{\mu}_{\text{aux}}^{\tilde{\sigma}} (dz) \\
& = \frac{1}{2} \left( \var_{\bar{\mu}_{\text{aux}}^{\tilde{\sigma}}} \left( y_s \right) + \var_{\bar{\mu}_{\text{aux}}^{\tilde{\sigma}}} \left( y_r \right) \right) \\
& = \frac{1}{2} \left( \var_{\mu_{\text{aux}}^{\tilde{\sigma}}} (y_s) +\var_{\mu_{\text{aux}}^{\tilde{\sigma}}} (y_r) \right) \overset{Lemma~\ref{l_variance_estimate}}{\lesssim} 1. \label{e_d_dt_int_aux3}
\end{align}
We also note that finite range interaction with strictly diagonal dominant condition
\begin{align}
\sum_{ 1 \leq |j-i| \leq R} |M_{ij}| + \delta \leq M_{ii} = 1 \qquad \text{for all } i \in [N]
\end{align}
imply there is a constant~$C$ such that for all~$i, j \in [N]$
\begin{align}
|M_{ij}| \leq \exp \left(-C|i-j|\right). \label{e_d_dt_int_aux4}
\end{align}
Plugging the estimates~\eqref{e_d_dt_int_aux1},~\eqref{e_d_dt_int_aux2},~\eqref{e_d_dt_int_aux3} and~\eqref{e_d_dt_int_aux4} into~\eqref{e_double_variable_covariance} yields
\begin{align}
\left| T_{\eqref{e_double_variable_covariance}} \right| &\lesssim \left( \sum_{\substack{ p \in \Lambda_1, \\ |p-i| \leq R }} \sum_{\substack{ r_1 \in E_n \\ s_1 \notin E_n \\ |r_1-s_1| \leq R  }} |M_{ip}| |M_{r_1 s_1}|  \exp\left(-C | p-r_1 | \right) \right) \\
& \quad \times \left( \sum_{\substack{ q \in \Lambda_1, \\ |q-j| \leq R }} \sum_{\substack{ r_2 \in E_n \\ s_2 \notin E_n \\ |r_2-s_2| \leq R  }} |M_{jq}| |M_{r_2 s_2}|  \exp\left(-C | q-r_2 | \right)\right) \\
& \lesssim \left(\sum_{\substack{ p \in \Lambda_1, \\ |p-i| \leq R }} \sum_{\substack{ r_1 \in E_n \\ s_1 \notin E_n \\ |r_1-s_1| \leq R  }} \exp\left(-C|i-p| \right)\exp\left(-C|r_1-s_1| \right)\exp\left(-C|p-r_1| \right) \right) \\
& \quad \times \left(\sum_{\substack{ q \in \Lambda_1, \\ |q-j| \leq R }} \sum_{\substack{ r_2 \in E_n \\ s_2 \notin E_n \\ |r_2-s_2| \leq R  }} \exp\left(-C|j-q| \right)\exp\left(-C|r_2-s_2| \right)\exp\left(-C|q-r_2| \right) \right) \\
& \lesssim R^2(R+2S)^2\exp\left(-CS\right) \lesssim S^2 \exp\left(-CS\right).
\end{align}
as desired. \qed

\medskip

It remains to provide the proof of Lemma~\ref{l_aux_cramers_representation}. This follows from a straightforward calculations. \\

\noindent \emph{Proof of Lemma~\ref{l_aux_cramers_representation}.} \ This follows from a direct computation. Indeed, it holds that
\begin{align}
&\left( \ln \int_{\frac{1}{K} \sum_{k \in \Lambda_1} x_k = \widetilde{m}} \exp \left( - H_{\text{aux}} \left(x^{\Lambda_1}, x^{\Lambda_2}\right) \right) \mathcal{L}(dx^{\Lambda_1})  \right.    \\
&\left. \qquad \quad + K \tilde{\sigma}\widetilde{m}- \ln \int_{\mathbb{R}^{\Lambda_1}} \exp \left(\tilde{\sigma}\sum_{k \in \Lambda_1} x_k - H_{\text{aux}} \left(x^{\Lambda_1}, x^{\Lambda_2}\right) \right) dx^{\Lambda_1}\right)\\
& \qquad = \ln \int_{\frac{1}{K} \sum_{k \in \Lambda_1} x_k = \widetilde{m}} \exp \left( K \tilde{\sigma}\widetilde{m} - H_{\text{aux}} \left(x^{\Lambda_1}, x^{\Lambda_2}\right) \right) \mathcal{L}(dx^{\Lambda_1}) \\
& \qquad \quad - \ln \int_{\mathbb{R}^{\Lambda_1}} \exp \left(\tilde{\sigma}\sum_{k \in \Lambda_1} x_k - H_{\text{aux}} \left(x^{\Lambda_1}, x^{\Lambda_2}\right) \right) dx^{\Lambda_1} \\
& \qquad = \ln \frac{\int_{\frac{1}{K} \sum_{k \in \Lambda_1} x_k = \widetilde{m}} \exp \left( \tilde{\sigma} \sum_{k \in \Lambda_1} x_k - H_{\text{aux}} \left(x^{\Lambda_1}, x^{\Lambda_2}\right) \right) \mathcal{L}(dx^{\Lambda_1})}{\int_{\mathbb{R}^{\Lambda_1}} \exp \left(\tilde{\sigma}\sum_{k \in \Lambda_1} x_k - H_{\text{aux}} \left(x^{\Lambda_1}, x^{\Lambda_2}\right) \right) dx^{\Lambda_1}} \\
& \qquad \overset{\eqref{e_m_and_mk}}{=} \ln \frac{\int_{\frac{1}{\sqrt{K}} \sum_{k \in \Lambda_1} \left(x_k - \widetilde{m}_k \right) =0} \exp \left( \tilde{\sigma} \sum_{k \in \Lambda_1} x_k - H_{\text{aux}} \left(x^{\Lambda_1}, x^{\Lambda_2}\right) \right) \mathcal{L}(dx^{\Lambda_1})}{\int_{\mathbb{R}^{\Lambda_1}} \exp \left(\tilde{\sigma}\sum_{k \in \Lambda_1} x_k - H_{\text{aux}} \left(x^{\Lambda_1}, x^{\Lambda_2}\right) \right) dx^{\Lambda_1}} = \ln g_{\text{aux}}(0).
\end{align}
\qed

\subsection{Proof of Lemma~\ref{l_marginal_otto_reznikoff_argument}} \label{s_marginal_otto}

Let us begin with computing the second derivatives of~$\bar{H}_{\Lambda_2}$. 

\begin{lemma} \label{l_second_der_bar_H}
Assume~$i \in \Lambda_2^{(n)}$ and~$j \in \Lambda_2^{(l)}$ with~$n \neq l$. Then it holds that
\begin{align}
\frac{d^2}{dx_i dx_j} \bar{H}_{\Lambda_2}\left(x^{\Lambda_2 }  \right) &= - \cov_{\mu_m \left(dx^{\Lambda_1} \big| x^{\Lambda_2}\right)} \left( \frac{\partial}{\partial x_i} H\left( x^{\Lambda_1}, x^{\Lambda_2} \right) -\frac{\partial}{\partial x_{i-R}} H\left( x^{\Lambda_1}, x^{\Lambda_2} \right) , \right. \\
& \qquad \qquad \qquad \qquad \qquad \left. \frac{\partial}{\partial x_j} H\left( x^{\Lambda_1}, x^{\Lambda_2} \right) -\frac{\partial}{\partial x_{j-R}} H\left( x^{\Lambda_1}, x^{\Lambda_2} \right)\right).
\end{align}
\end{lemma}

\noindent \emph{Proof of Lemma~\ref{l_second_der_bar_H}.} \  
Recall the definition~\eqref{d_def_marginal_hamiltonian1} of~$\bar{H}_{\Lambda_2}\left(x^{\Lambda_2 }  \right)$ given by
\begin{align}
\bar{H}_{\Lambda_2}\left(x^{\Lambda_2 }  \right) &= -\ln \int_{\frac{1}{N} \sum_{k \in \Lambda_1} x_k = m - \frac{1}{N} \sum_{i \in \Lambda_2} x_i} \exp\left( - H\left(x^{\Lambda_1}, x^{\Lambda_2}\right) \right)  \mathcal{L}^{K-1}(dx^{\Lambda_1}).
\end{align}
Fix~$x^{\Lambda_2} \in \mathbb{R}^{\Lambda_2}$ and define a vector~$z ^{\Lambda_1} \in \mathbb{R}^{\Lambda_1}$ as
\begin{align} \label{e_def_linear_trans_z}
z_k : = \begin{cases} x_{k+R} \qquad &\text{if } k+R \in \Lambda_2 ,\\
0 \qquad &\text{otherwise},
\end{cases} \qquad \text{for each } k \in \Lambda_1,
\end{align}
and let
\begin{align} \label{e_def_linear_trans_y}
y^{\Lambda_1} : = x^{\Lambda_1} + z^{\Lambda_1}.
\end{align}
Note that~$\sum_{k \in \Lambda_1} z_k = \sum_{i \in \Lambda_2} x_i$. In particular, it holds that
\begin{align}
\sum_{k=1}^N x_k = \sum_{k \in \Lambda_1} x_k + \sum_{i \in \Lambda_2} x_i = \sum_{k \in \Lambda_1} x_k + \sum_{k \in \Lambda_1} z_k = \sum_{k \in \Lambda_1 }y_k. \label{e_change_of_var_aux}
\end{align}
With this observation, it follows from change of variables that
\begin{align}
\bar{H}_{\Lambda_2}\left(x^{\Lambda_2 }  \right) & = - \ln \int_{\frac{1}{N} \sum_{k \in \Lambda_1} y_k =m } \exp\left(- H\left( y^{\Lambda_1} - z^{\Lambda_1}, x^{\Lambda_2} \right) \right) \mathcal{L}^{K-1}(dy^{\Lambda_1}). \label{e_bar_h_lambda_2}
\end{align}
Note also that a direct calculation yields for~$i \in \Lambda_2$,
\begin{align}
\frac{d}{dx_i}  H \left( y^{\Lambda_1} - z^{\Lambda_1}, x^{\Lambda_2} \right)  & =  \frac{\partial}{\partial x_i}H \left( y^{\Lambda_1} - z^{\Lambda_1}, x^{\Lambda_2} \right) - \sum_{k \in \Lambda_1}  \frac{\partial}{\partial x_k} H \left( y^{\Lambda_1} - z^{\Lambda_1}, x^{\Lambda_2} \right) \cdot \frac{dz_k}{dx_i}  \\
& = \frac{\partial}{\partial x_i}H \left( y^{\Lambda_1} - z^{\Lambda_1}, x^{\Lambda_2} \right) - \frac{\partial}{\partial x_{i-R}}H \left( y^{\Lambda_1} - z^{\Lambda_1}, x^{\Lambda_2} \right) \\
& = \frac{\partial}{\partial x_i} H\left( x^{\Lambda_1}, x^{\Lambda_2} \right) -\frac{\partial}{\partial x_{i-R}} H\left( x^{\Lambda_1}, x^{\Lambda_2} \right) .\label{e_derivative_linear_transformation}
\end{align}
Then a combination of~\eqref{e_change_of_var_aux},~\eqref{e_bar_h_lambda_2} and~\eqref{e_derivative_linear_transformation} followed by change of variables yields
\begin{align}
&\frac{d}{dx_i} \bar{H}_{\Lambda_2}\left(x^{\Lambda_2 }  \right)\\
& = - \frac{\int_{\frac{1}{N} \sum_{k \in \Lambda_1} y_k =m } \frac{d}{dx_i}\left(- H\left( y^{\Lambda_1} - z^{\Lambda_1}, x^{\Lambda_2} \right) \right)\exp\left(- H\left( y^{\Lambda_1} - z^{\Lambda_1}, x^{\Lambda_2} \right) \right) \mathcal{L}^{K-1}(dy^{\Lambda_1})}{\int_{\frac{1}{N} \sum_{k \in \Lambda_1} y_k =m } \exp\left(- H\left( y^{\Lambda_1} - z^{\Lambda_1}, x^{\Lambda_2} \right) \right) \mathcal{L}^{K-1}(dy^{\Lambda_1})}  \\
& = - \frac{\int_{\frac{1}{N} \sum_{k =1}^{N} x_k =m } \left(-\frac{\partial}{\partial x_i} H\left( x^{\Lambda_1}, x^{\Lambda_2} \right) +\frac{\partial}{\partial x_{i-R}} H\left( x^{\Lambda_1}, x^{\Lambda_2} \right) \right)\exp\left( - H\left(x^{\Lambda_1}, x^{\Lambda_2}\right) \right)  \mathcal{L}^{K-1}(dx^{\Lambda_1})}{\int_{\frac{1}{N} \sum_{k =1}^{N} x_k = m } \exp\left( - H\left(x^{\Lambda_1}, x^{\Lambda_2}\right) \right)  \mathcal{L}^{K-1}(dx^{\Lambda_1})} \\
& \overset{\eqref{e_def_tilde_m}}{=} - \frac{\int_{\frac{1}{K} \sum_{k \in \Lambda_1} x_k =\widetilde{m} } \left(-\frac{\partial}{\partial x_i} H\left( x^{\Lambda_1}, x^{\Lambda_2} \right) +\frac{\partial}{\partial x_{i-R}} H\left( x^{\Lambda_1}, x^{\Lambda_2} \right) \right)\exp\left( - H\left(x^{\Lambda_1}, x^{\Lambda_2}\right) \right)  \mathcal{L}^{K-1}(dx^{\Lambda_1})}{\int_{\frac{1}{N} \sum_{k =1}^{N} x_k = m } \exp\left( - H\left(x^{\Lambda_1}, x^{\Lambda_2}\right) \right)  \mathcal{L}^{K-1}(dx^{\Lambda_1})} \\
& = \mathbb{E}_{\mu_m \left(dx^{\Lambda_1} \big| x^{\Lambda_2}\right)} \left[ \frac{\partial}{\partial x_i} H\left( x^{\Lambda_1}, x^{\Lambda_2} \right) -\frac{\partial}{\partial x_{i-R}} H\left( x^{\Lambda_1}, x^{\Lambda_2} \right) \right].
\end{align}
Let us turn to the computation of the second derivative of~$\bar{H}$. Observe that for any distinct~$k, l \in [N]$,
\begin{align}
\frac{\partial^2}{\partial x_k \partial x_l} H\left(x^{\Lambda_1}, x^{\Lambda_2} \right) = M_{kl}.
\end{align}
In particular if~$|k-l|>R$, we have
\begin{align}
\frac{\partial^2}{\partial x_k \partial x_l} H\left(x^{\Lambda_1}, x^{\Lambda_2} \right) = M_{kl} = 0.
\end{align}
Suppose~$i \in \Lambda_2^{(n)}$ and~$j \in \Lambda_2^{(l)}$ with~$n \neq l$. For~$L> 2R$ (see Figure~\ref{f_two_scale_decomposition_of_lattice}), it holds that
$|i-j| \geq L >2R$. This implies
\begin{align}
\text{dist} \left( \{i, i-R\}, \{j, j-R\}  \right) >R.
\end{align}
With this observation, a similar computation from above yields the desired formula
\begin{align}
\frac{d^2}{dx_i dx_j} \bar{H}_{\Lambda_2}\left(x^{\Lambda_2 }  \right) 
& = - \cov_{\mu_m \left(dx^{\Lambda_1} \big| x^{\Lambda_2}\right)} \left( \frac{\partial}{\partial x_i} H\left( x^{\Lambda_1}, x^{\Lambda_2} \right) -\frac{\partial}{\partial x_{i-R}} H\left( x^{\Lambda_1}, x^{\Lambda_2} \right)  ,\right. \\
& \qquad \qquad \qquad \qquad \qquad \left. \frac{\partial}{\partial x_j} H\left( x^{\Lambda_1}, x^{\Lambda_2} \right) -\frac{\partial}{\partial x_{j-R}} H\left( x^{\Lambda_1}, x^{\Lambda_2} \right)\right).
\end{align}
\qed 

\medskip

The main ingredient for proving Lemma~\ref{l_marginal_otto_reznikoff_argument} is the decay of correlations. We consider the conditional measure~$\mu_m \left( dx^{\Lambda_1} \big| x^{\Lambda_2} \right)$ which is a ce with mean spin~$\widetilde{m}$ given by~\eqref{e_def_tilde_m}. We define the corresponding conditional gce~$\mu ^{\tau} \left(dx^{\Lambda_1} \big| x^{\Lambda_2} \right)$ as
\begin{align} \label{e_def_conditional_gce}
\mu^{\tau}\left(dx^{\Lambda_1} \big| x^{\Lambda_2} \right) : = \frac{1}{Z} \exp \left( \tau \sum_{k \in \Lambda_1} x_k - H\left(x^{\Lambda_1}, x^{\Lambda_2} \right)  \right) dx^{\Lambda_1},
\end{align}
where~$\tau = \tau(\widetilde{m})$ is given by (cf. Definition~\ref{a_relation_sigma_m})
\begin{align}
\widetilde{m} = \frac{1}{K} \int \sum_{k \in \Lambda_1} x_k \mu ^{\tau}\left(dx^{\Lambda_1} \big| x^{\Lambda_2} \right).
\end{align}
Note that~$\mu^{\tau}\left(dx^{\Lambda_1} \big| x^{\Lambda_2} \right)$ and~$\mu_{m}\left(dx^{\Lambda_1} \big| x^{\Lambda_2}\right)$ satisfy the same structural assumptions made in Proposition~\ref{l_decay_corr_ce}. As a consequence, Proposition~\ref{l_decay_corr_ce}, i.e. decay of correlations, also holds for $\mu^{\tau}\left(dx^{\Lambda_1} \big| x^{\Lambda_2} \right)$ and $\mu_{m}\left(dx^{\Lambda_1} \big| x^{\Lambda_2}\right)$. Now we are ready to give a proof of Lemma~\ref{l_marginal_otto_reznikoff_argument}. \\

\noindent \emph{Proof of Lemma~\ref{l_marginal_otto_reznikoff_argument}.} \  Assume~$ n \neq l$ and let~$i \in \Lambda_2^{(n)}$,~$j\in \Lambda_2^{(l)}$ be given. It is enough to prove that
\begin{align}
\left| \frac{d^2}{dx_i dx_j} \bar{H}_{\Lambda_2}\left(x^{\Lambda_2 }  \right) \right| \lesssim \frac{1}{K^{1-\varepsilon}} +  \exp \left( -CL|n-l| \right).
\end{align}
Note that for each~$k \in [N]$
\begin{align}
\frac{\partial}{\partial x_k} H \left(x^{\Lambda_1}, x^{\Lambda_2} \right) = x_k + \psi_b ' (x_k) + s_k + \sum_{l : \ 1\leq |l-k| \leq R} M_{kl}x_l.
\end{align}
In particular for~$L$ large enough we have
\begin{align}
&\text{dist} \left( \supp \left( \frac{\partial}{\partial x_i} H\left( x^{\Lambda_1}, x^{\Lambda_2} \right) -\frac{\partial}{\partial x_{i-R}} H\left( x^{\Lambda_1}, x^{\Lambda_2} \right) \right) , \right. \\
&\left.  \qquad \quad  \supp \left( \frac{\partial}{\partial x_j} H\left( x^{\Lambda_1}, x^{\Lambda_2} \right) -\frac{\partial}{\partial x_{j-R}} H\left( x^{\Lambda_1}, x^{\Lambda_2} \right) \right) \right) \geq (L-2R) |n-l| \geq \frac{L}{2} |n-l|. \label{e_distance_supp}
\end{align}
Recall the definition~\eqref{e_def_conditional_gce} of~$\mu^{\tau}$. Because it is the gce on the one-dimensional lattice, one has the uniform moment estimate as in Lemma~\ref{l_moment_estimate}, i.e. for each~$n \in \mathbb{N}$, there is a constant~$C=C(n)$ such that for each $i \in \Lambda_1$,
\begin{align}
\mathbb{E}_{\mu^{\tau}} \left[ \left| X_i - \mathbb{E}_{\mu^{\tau}}[X_i] \right|^n \right] \leq C(n) .\label{e_conditional_moment_estimate}
\end{align}

Therefore a combination of Proposition~\ref{l_decay_corr_ce}, Lemma~\ref{l_second_der_bar_H}, the fact that covariances are invariant under the addition of constants,~\eqref{e_distance_supp}, and~\eqref{e_conditional_moment_estimate} yields the desired estimate
\begin{align}
\left| \frac{d^2}{dx_i dx_j} \bar{H}_{\Lambda_2}\left(x^{\Lambda_2 }  \right) \right| \lesssim \frac{1}{K^{1-\varepsilon}} +  \exp \left( -CL|n-l| \right).
\end{align}
\qed

\medskip

\section{Proof of Proposition~\ref{p_combined_lsi} } \label{s_proof_combined_lsi}
The proof is a slight adaptation of the argument for the two-scale criterion for LSI (cf. Theorem~\ref{a_two_scale},~\cite[Theorem 3]{GrOtViWe09}). For the convenience of the reader we give all details. Let~$\Phi : \mathbb{R} \to \mathbb{R}$ denote the function~$\Phi (x) = x \ln x$. Observe that we can decompose the relative entropy into
\begin{align}
&\int \Phi (f(x)) \mu_m (dx) - \Phi \left( \int f(x) \mu_m (dx) \right) \\
&\qquad = \int \left[ \int \Phi \left(f(x)\right)\mu_m \left(dx^{\Lambda_1} \big| x^{\Lambda_2} \right) - \Phi \left(\int f(x) \mu_m \left(dx^{\Lambda_1} \big| x^{\Lambda_2} \right) \right) \right] \bar{\mu}_m \left(dx^{\Lambda_2}\right) \label{e_combined_lsi_1}\\
& \qquad \quad  + \int \Phi \left( \bar{f}\left(x^{ \Lambda_2} \right)   \right) \bar{\mu}_m \left(dx^{\Lambda_2} \right) - \Phi \left( \int \bar{f} \left( x^{\Lambda_2} \right) \bar{\mu}_m \left( dx^{\Lambda_2} \right) \right), \label{e_combined_lsi_2}
\end{align}
where~$\bar{f}\left( x^{\Lambda_2} \right) : = \int f(x) \mu_m \left( dx^{\Lambda_1} \big| x^{\Lambda_2} \right)$.
Let us begin with estimation of~\eqref{e_combined_lsi_1}. As~$\mu_m \left( dx^{\Lambda_1} \big| x^{\Lambda_2} \right)$ satisfies LSI$(\rho_1)$, it follows that
\begin{align}
&\int \left[ \int \Phi \left(f(x)\right)\mu_m \left(dx^{\Lambda_1} \big| x^{\Lambda_2} \right) - \Phi \left(\int f(x) \mu_m \left(dx^{\Lambda_1} \big| x^{\Lambda_2} \right) \right) \right] \bar{\mu}_m \left(dx^{\Lambda_2}\right) \\
& \qquad \leq \int \left[  \frac{1}{2\rho_1} \int \frac{ \left| \nabla_{\Lambda_1} f \right|^2}{f} \mu_m \left(dx^{\Lambda_1} \big| x^{\Lambda_2} \right) \right] \bar{\mu}_m \left(dx^{\Lambda_2} \right) \\
& \qquad = \frac{1}{2\rho_1} \int \frac{\left| \nabla_{\Lambda_1} f \right|^2 }{f} \mu_m \left(dx \right).
\end{align}
Let us turn to the estimation of~\eqref{e_combined_lsi_2}. It holds by the assumption on $\bar{\mu}_m \left( dx^{\Lambda_2} \right)$ satisfying LSI$(\rho_2)$ that
\begin{align}
\int \Phi \left( \bar{f}\left(x^{ \Lambda_2} \right)   \right) \bar{\mu}_m \left(dx^{\Lambda_2} \right) - \Phi \left( \int \bar{f} \left( x^{\Lambda_2} \right) \bar{\mu}_m \left( dx^{\Lambda_2} \right) \right) \leq \frac{1}{2 \rho_2} \int \frac{\left| \nabla_{\Lambda_2} \bar{f}\right|^2 }{\bar{f} }\bar{\mu}_m \left(dx^{\Lambda_2} \right). \label{e_estimation_last_term}
\end{align}

\begin{lemma}[Analogue of Lemma 21 in~\cite{GrOtViWe09}] \label{l_derivative_of_bar_f}
It holds that
\begin{align}
 \nabla_{\Lambda_2} \bar{f}  &=\left( \mathbb{E}_{\mu_m \left(dx^{\Lambda_1} \big| x^{\Lambda_2}\right)} \left[ \frac{\partial}{\partial x_i} f(x) - \frac{\partial}{\partial x_{i-R}} f(x)   \right] \right)_{i \in \Lambda_2} \\
& \quad - \cov_{\mu_m \left(dx^{\Lambda_1} \big| x^{\Lambda_2} \right)}\left( f(x) , \left(\frac{\partial}{\partial x_i} H\left(x^{\Lambda_1}, x^{\Lambda_2}\right) -\frac{\partial}{\partial x_{i-R}} H\left(x^{\Lambda_1}, x^{\Lambda_2}\right) \right)_{i \in \Lambda_2} \right).
\end{align}
\end{lemma}

\noindent \emph{Proof of Lemma~\ref{l_derivative_of_bar_f}.} \  Given~$x^{\Lambda_2} \in \mathbb{R}^{\Lambda_2}$, recall the definition~\eqref{e_def_linear_trans_z} and~\eqref{e_def_linear_trans_y} of~$z^{\Lambda_1}$ and~$y^{\Lambda_1}$, respectively. Similar calculations as in Lemma~\ref{l_second_der_bar_H} implies
\begin{align}
\frac{\partial}{\partial x_i} \bar{f} \left(x^{\Lambda_2} \right) &= \frac{\partial}{\partial x_i} \frac{\int_{ \frac{1}{N}\sum_{k \in \Lambda_1} y_k = m} f\left( y^{\Lambda_1}-z^{\Lambda_1} , x^{\Lambda_2} \right) \exp \left( - H\left(y^{\Lambda_1}-z^{\Lambda_1} , x^{\Lambda_2} \right)\right) \mathcal{L}(dy ^{\Lambda_{1}}) }{\int_{ \frac{1}{N}\sum_{k \in \Lambda_1} y_k = m} \exp \left( - H\left(y^{\Lambda_1}-z^{\Lambda_1} , x^{\Lambda_2} \right)\right) \mathcal{L}(dy ^{\Lambda_{1}})} \\
& = \mathbb{E}_{\mu_m \left(dx^{\Lambda_1} \big| x^{\Lambda_2}\right)} \left[ \frac{\partial}{\partial x_i} f(x) -\frac{\partial}{\partial x_{i-R}} f(x)   \right]  \\
&\quad - \cov_{\mu_m \left(dx^{\Lambda_1} \big| x^{\Lambda_2} \right)}\left( f(x) , \frac{\partial}{\partial x_i} H\left(x^{\Lambda_1}, x^{\Lambda_2}\right) -\frac{\partial}{\partial x_{i -R}} H\left(x^{\Lambda_1}, x^{\Lambda_2}\right)\right).
\end{align}
\qed

\medskip

The next ingredient is the estimation of covariances.

\begin{lemma}[Lemma 22 in~\cite{GrOtViWe09}]\label{l_covariance_lsi_estimate}
Let~$\mu$ be a probability measure on an Euclidean space~$X$. Suppose~$\mu$ satisfies LSI$(\rho)$ for some~$\rho>0$. Then for any two Lipschitz functions~$f : X \to [0, +\infty)$ and~$g : X \to \mathbb{R}$,
\begin{align}
\left| \cov_{\mu} \left( f, g\right) \right| \leq \frac{\| \nabla g \|_{L^{\infty}(\mu)}}{\rho} \sqrt{ \left( \int f d\mu \right) \left( \int \frac{|\nabla f |^2}{f} d\mu \right)}.
\end{align}
\end{lemma}

Now we estimate the integrand in the right hand side of~\eqref{e_estimation_last_term} with the help of Lemma~\ref{l_derivative_of_bar_f} and Lemma~\ref{l_covariance_lsi_estimate}.

\begin{lemma}\label{l_Jensen}
It holds that
\begin{align}
\frac{\left| \nabla_{\Lambda_2} \bar{f} \right|^2}{\bar{f}} \lesssim \frac{1}{\rho_1 ^2}  \int \frac{ \left| \nabla_{\Lambda_1} f \right|^2}{f} \mu_m \left( dx^{\Lambda_1} \big| x^{\Lambda_2} \right)+ \int \frac{ \left| \nabla_{\Lambda_2} f \right|^2}{f} \mu_m \left( dx^{\Lambda_1} \big| x^{\Lambda_2} \right) .
\end{align}
\end{lemma}

\noindent \emph{Proof of Lemma~\ref{l_Jensen}.} \  
An application of Young's inequality yields
\begin{align}
\frac{\left| \nabla_{\Lambda_2} \bar{f} \right|^2}{\bar{f}} & \overset{Lemma~\ref{l_derivative_of_bar_f} }{=} \frac{1}{\bar{f}} \left| \left( \mathbb{E}_{\mu_m \left(dx^{\Lambda_1} \big| x^{\Lambda_2}\right)} \left[ \frac{\partial}{\partial x_i} f(x) - \frac{\partial}{\partial x_{i-R}} f(x)   \right] \right)_{i \in \Lambda_2}\right. \\
& \qquad \quad  \quad \left. - \cov_{\mu_m \left(dx^{\Lambda_1} \big| x^{\Lambda_2} \right)}\left( f(x) , \left(\frac{\partial}{\partial x_i} H\left(x^{\Lambda_1}, x^{\Lambda_2}\right) -\frac{\partial}{\partial x_{i-R}} H\left(x^{\Lambda_1}, x^{\Lambda_2}\right) \right) \right)_{i \in \Lambda_2} \right|^2 \\
&  \quad  \leq \frac{4}{\bar{f}} \left| \left( \mathbb{E}_{\mu_m \left(dx^{\Lambda_1} \big| x^{\Lambda_2}\right)} \left[ \frac{\partial}{\partial x_i} f(x) \right] \right)_{i \in \Lambda_2} \right|^2  \label{e_young1} \\
& \quad \quad + \frac{4}{\bar{f}} \left| \left( \mathbb{E}_{\mu_m \left(dx^{\Lambda_1} \big| x^{\Lambda_2}\right)} \left[ \frac{\partial}{\partial x_{i-R}} f(x) \right] \right)_{i \in \Lambda_2} \right|^2 \label{e_young2} \\
& \quad \quad + \frac{4}{\bar{f}} \left| \cov_{\mu_m \left(dx^{\Lambda_1} \big| x^{\Lambda_2} \right)}\left( f(x) , \frac{\partial}{\partial x_i} H\left(x^{\Lambda_1}, x^{\Lambda_2}\right)  \right)_{i \in \Lambda_2} \right|^2 \label{e_young3} \\
& \quad \quad + \frac{4}{\bar{f}} \left| \cov_{\mu_m \left(dx^{\Lambda_1} \big| x^{\Lambda_2} \right)}\left( f(x) , \frac{\partial}{\partial x_{i-R}} H\left(x^{\Lambda_1}, x^{\Lambda_2}\right)  \right)_{i \in \Lambda_2} \right|^2. \label{e_young4}
\end{align}
Let us begin with estimation of~\eqref{e_young1}. Cauchy's inequality implies that
\begin{align}
\left|  \mathbb{E}_{\mu_m \left(dx^{\Lambda_1} \big| x^{\Lambda_2}\right)} \left[ \nabla_{\Lambda_2} f(x) \right] \right|^2 &\leq \int f(x) \mu_m \left( dx^{\Lambda_1} \big| x^{\Lambda_2} \right) \int \frac{ \left| \nabla_{\Lambda_2} f \right|^2}{f} \mu_m \left( dx^{\Lambda_1} \big| x^{\Lambda_2} \right) \\
& = \bar{f}\left(x^{\Lambda_2} \right) \int \frac{ \left| \nabla_{\Lambda_2} f \right|^2}{f} \mu_m \left( dx^{\Lambda_1} \big| x^{\Lambda_2} \right), 
\end{align}
and as a consequence
\begin{align}
T_{\eqref{e_young1}} \leq 4 \int \frac{ \left| \nabla_{\Lambda_2} f \right|^2}{f} \mu_m \left( dx^{\Lambda_1} \big| x^{\Lambda_2} \right).
\end{align}
Let us turn to the estimation of~\eqref{e_young2}. Note that~$i-R \in \Lambda_1$ for each~$i \in \Lambda_2$. Then a similar computation from above yields
\begin{align}
T_{\eqref{e_young2}} \leq \frac{4}{\bar{f}} \left|\mathbb{E}_{\mu_m \left(dx^{\Lambda_1} \big| x^{\Lambda_2}\right)} \left[ \nabla_{\Lambda_1} f(x) \right] \right|^2 \leq 4 \int \frac{ \left| \nabla_{\Lambda_1} f \right|^2}{f} \mu_m \left( dx^{\Lambda_1} \big| x^{\Lambda_2} \right).
\end{align}
Let us turn to the estimation of~\eqref{e_young3}.
\begin{align}
&\left| \cov_{\mu_m \left(dx^{\Lambda_1} \big| x^{\Lambda_2} \right)}\left( f(x) , \frac{\partial}{\partial x_i} H\left(x^{\Lambda_1}, x^{\Lambda_2}\right)  \right)_{i \in \Lambda_2} \right|^2 \\
& \qquad = \sup_{\substack{ |y|\leq 1 \\ y \in \mathbb{R}^{\Lambda_2} }}\left[ \cov_{\mu_m \left(dx^{\Lambda_1} \big| x^{\Lambda_2} \right)}\left( f(x) ,  \nabla_{\Lambda_2} H\left(x^{\Lambda_1}, x^{\Lambda_2}\right) \cdot y \right) \right]^2 \\
& \overset{Proposition~\ref{p_conditional_lsi}, Lemma~\ref{l_covariance_lsi_estimate}}{\leq} \frac{1}{\rho_1 ^2} \left( \sup_{\substack{ |y|\leq 1 \\ y \in \mathbb{R}^{\Lambda_2} }} \sup_{\substack{ x^{\Lambda_1} \in \mathbb{R}^{\Lambda_1} }} \left| \nabla_{\Lambda_1} \left( \nabla_{\Lambda_2} H\left(x^{\Lambda_1}, x^{\Lambda_2}\right) \cdot y\right)  \right|^2  \right) \\
& \qquad \qquad \qquad \times \left( \int f(x) \mu_m \left( dx^{\Lambda_1} \big| x^{\Lambda_2} \right) \right) \left( \int \frac{ \left| \nabla_{\Lambda_1} f \right|^2}{f} \mu_m \left( dx^{\Lambda_1} \big| x^{\Lambda_2} \right)\right) \\
& \qquad \lesssim \frac{1}{\rho_1 ^2} \bar{f}\left(x^{\Lambda_2} \right)\left( \int \frac{ \left| \nabla_{\Lambda_1} f \right|^2}{f} \mu_m \left( dx^{\Lambda_1} \big| x^{\Lambda_2} \right)\right)  .
\end{align}
Therefore
\begin{align}
T_{\eqref{e_young3}} \lesssim \frac{1}{\rho_1 ^2}  \int \frac{ \left| \nabla_{\Lambda_1} f \right|^2}{f} \mu_m \left( dx^{\Lambda_1} \big| x^{\Lambda_2} \right) .
\end{align}
Similar computation also yields
\begin{align}
T_{\eqref{e_young4}} \lesssim \frac{1}{\rho_1 ^2}  \int \frac{ \left| \nabla_{\Lambda_1} f \right|^2}{f} \mu_m \left( dx^{\Lambda_1} \big| x^{\Lambda_2} \right) .
\end{align}
\qed 

\medskip

Let us now give a proof of Proposition~\ref{p_combined_lsi}. \\ 

\noindent \emph{Proof of Proposition~\ref{p_combined_lsi}.} \ A combination of~\eqref{e_estimation_last_term} and Lemma~\ref{l_Jensen} yields that
\begin{align}
T_{\eqref{e_combined_lsi_2}} &\lesssim \frac{1}{2 \rho_2} \max \left \{ 1, \frac{1}{\rho_1 ^2} \right \} \int \left( \int \frac{\left| \nabla f \right|^2}{f} \mu_m \left(dx^{\Lambda_1} \big| x^{\Lambda_2} \right) \right) \bar{\mu}_m \left(dx^{\Lambda_2} \right)  \\
& = \frac{1}{2 \rho_2} \max \left \{ 1, \frac{1}{\rho_1 ^2} \right \} \int \frac{\left| \nabla f \right|^2}{f} \mu_m \left(dx\right).
\end{align}
Therefore
\begin{align}
\int \Phi (f) d\mu - \Phi \left( \int f d\mu \right) &= T_{\eqref{e_combined_lsi_1}} + T_{\eqref{e_combined_lsi_2}} \\
& \lesssim\frac{1}{2\rho_1} \int \frac{ \left|\nabla_{\Lambda_1}f \right|^2}{f} \mu_m \left(dx \right)+ \frac{1}{2 \rho_2} \max \left \{ 1, \frac{1}{\rho_1 ^2} \right \} \int \frac{\left| \nabla f \right|^2}{f} \mu_m \left(dx\right) \\
& \leq \frac{1}{2}\left(\frac{1}{\rho_1} + \frac{1}{\rho_2}\max \left \{ 1, \frac{1}{\rho_1 ^2} \right \} \right) \int \frac{\left| \nabla f \right|^2}{f} \mu_m \left(dx\right).
\end{align}
This finishes the proof of Proposition~\ref{p_combined_lsi}.
\qed 

\medskip

\appendix                                     
\section{Criteria for the logarithmic Sobolev inequality} \label{a_lsi_criteria}
In this section we state several standard criteria for deducing a LSI. For proofs we refer to the literature. For a general introduction and more comments on the LSI we refer the reader to~\cite{Led01,Led01a,Roy99,BaGeLe14}. 

\begin{theorem}[Tensorization Principle~\cite{Gro75}] \label{a_tensorization}
Let~$\mu_1$ and~$\mu_2$ be probability measures on Euclidean spaces~$X_1$ and~$X_2$ respectively. Suppose that~$\mu_1$ and~$\mu_2$ satisfy LSI$(\rho_1)$ and LSI$(\rho_2)$ respectively. Then the product measure~$\mu_1 \bigotimes \mu_2$ satisfies LSI$(\rho)$, where~$\rho = \min \{\rho_1, \rho_2\}$.
\end{theorem}

\begin{theorem}[Holley-Stroock Perturbation Principle~\cite{HolStr87}] \label{a_holley_stroock}
Let~$\mu_1$ be a probability measure on Euclidean space~$X$ and~$\delta \psi : X \to \mathbb{R}$ be a bounded function. Define a probability measure~$\mu_2$ on~$X$ by
\begin{align}
\mu_2 (dx) : = \frac{1}{Z} \exp \left( - \delta \psi (x) \right) \mu_1 (dx).
\end{align}
Suppose that~$\mu_1$ satisfies LSI$(\rho_1)$. Then~$\mu_2$ also satisfies LSI with constant
\begin{align}
\rho_2 = \rho_1 \exp \left( - \text{osc } \delta \psi \right),
\end{align}
where~$\text{osc } \delta \psi : = \sup \delta \psi - \inf \delta \psi $.
\end{theorem}
\begin{theorem}[Bakry-\'{E}mery criterion~\cite{BaEm85}]\label{a_bakry_emery} Let~$X$ be a~$N$-dimensional Euclidean space and $H \in C^2 (X)$. Define a probability measure~$\mu$ on~$X$ by
\begin{align}
\mu(dx): = \frac{1}{Z} \exp\left( - H(x)\right) dx.
\end{align}
Suppose there is a constant~$\rho>0$ such that~$\Hess H \geq \rho$. More precisely, for all~$u, v \in X$,
\begin{align}
\langle v, \Hess H(u) v \rangle \geq \rho |v|^2.
\end{align}
Then~$\mu$ satisfies LSI$(\rho)$.
\end{theorem}
\begin{theorem}[Otto-Reznikoff Criterion~\cite{OttRez07}] \label{a_otto_reznikoff} Let~$X = X_1 \times \cdots \times X_N$ be a direct product of Euclidean spaces and~$H \in C^2 (X)$. Define a probability measure~$\mu$ on~$X$ by
\begin{align}
\mu(dx) : = \frac{1}{Z} \exp \left( - H(x) \right) dx.
\end{align}
Assume that
\begin{itemize}
\item For each~$i \in \{1, \cdots, N\}$, the conditional measures~$\mu (dx_i | \bar{x}_i )$ satisfy LSI$(\rho_i)$.
\item For each~$ 1 \leq i \neq j \leq N$ there is a constant~$\kappa_{ij} \in (0, \infty)$ with
\begin{align}
\left| \nabla_i \nabla_j H(x) \right| \leq \kappa_{ij}. \qquad \text{for all } x \in X.
\end{align}
Here,~$| \cdot |$ denotes the operator norm of a bilinear form.
\item Define a symmetric matrix~$A= (A_{ij})_{1 \leq i, j \leq N}$ by
\begin{align}
A_{ij} = \begin{cases} \rho_i, \qquad &\text{if } i=j \\
- \kappa_{ij}, \qquad &\text{if } i \neq j
\end{cases}.
\end{align}
Assume that there is a constant~$\rho \in (0, \infty)$ with
\begin{align}
A \geq \rho \Id,
\end{align}
in the sense of quadratic forms.
\end{itemize}
Then~$\mu$ satisfies LSI$(\rho)$.
\end{theorem}

\begin{theorem}[Two-Scale Criterion~\cite{GrOtViWe09}] \label{a_two_scale}
Let~$X$ and~$Y$ be Euclidean spaces. Consider a probability measure~$\mu$ on~$X$ defined by
\begin{align}
\mu(dx) : = \frac{1}{Z} \exp\left(-H(x)\right)dx.
\end{align}
Let~$P : X \to Y$ be a linear operator such that for some~$N \in \mathbb{N}$,
\begin{align}
PNP^t = \Id_Y.
\end{align}
Define
\begin{align}
\kappa : = \max \left \{ \langle \Hess H(x) \cdot u, v \rangle : \ u \in \text{Ran}(NP^t P), v \in \text{Ran}(\Id_X - NP^t P), |u|=|v|=1    \right \}.
\end{align}
Assume that
\begin{itemize}
\item $\kappa < \infty$ 
\item There is~$\rho_1 \in (0, \infty)$ such that the conditional measure~$\mu(dx | Px=y)$ satisfies LSI$(\rho_1)$ for all~$y \in Y$.
\item There is~$\rho_2 \in (0, \infty)$ such that the marginal measure~$\bar{\mu}= P_{\#}\mu$ satisfies LSI$(\rho_2 N)$.
\end{itemize}
Then~$\mu$ satisfies LSI$(\rho)$, where
\begin{align}
\rho : = \frac{1}{2} \left( \rho_1 +\rho_2 + \frac{ \kappa^2}{\rho_1} - \sqrt{ \left( \rho_1 + \rho_2 + \frac{\kappa ^2 }{\rho_1} \right)^2 - 4 \rho_1 \rho_2 }    \right) >0.
\end{align}
\end{theorem}

\section*{Acknowledgment}
This research has been partially supported by NSF grant DMS-1407558. The authors want to thank Felix Otto and H.T.~Yau for bringing this problem to their attention. The authors are also thankful to many people discussing the problem and helping to improve the preprint. Among them are Tim Austin, Frank Barthe, Marek Biskup, Pietro Caputo, Jean-Dominique Deuschel, Max Fathi, Andrew Krieger, Michel Ledoux, Sangchul Lee, Thomas Liggett, Felix Otto, Daniel Ueltschi, and Tianqi Wu. The authors want to thank Marek Biskup, UCLA and KFAS for financial support.

\bibliographystyle{alpha}
\bibliography{bib}

\begin{thebibliography}{GOVW09}

\bibitem[BE85]{BaEm85}
D.~Bakry and Michel \'Emery.
\newblock Diffusions hypercontractives.
\newblock In {\em S\'eminaire de probabilit\'es, {XIX}, 1983/84}, volume 1123
  of {\em Lecture Notes in Math.}, pages 177--206. Springer, Berlin, 1985.

\bibitem[BGL14]{BaGeLe14}
Dominique Bakry, Ivan Gentil, and Michel Ledoux.
\newblock {\em Analysis and geometry of {M}arkov diffusion operators}, volume
  348 of {\em Grundlehren der Mathematischen Wissenschaften [Fundamental
  Principles of Mathematical Sciences]}.
\newblock Springer, Cham, 2014.

\bibitem[BL76]{Bra76}
Herm~Jan Brascamp and Elliott~H Lieb.
\newblock Best constants in {Y}oung's inequality, its converse, and its
  generalization to more than three functions.
\newblock {\em Advances in Mathematics}, 20(2):151 -- 173, 1976.

\bibitem[Cap03]{Cap03}
Pietro Caputo.
\newblock Uniform {P}oincar\'e inequalities for unbounded conservative spin
  systems: the non-interacting case.
\newblock {\em Stochastic Process. Appl.}, 106(2):223--244, 2003.

\bibitem[CFMP05]{CaFeMePr05}
Marzio Cassandro, Pablo~Augusto Ferrari, Immacolata Merola, and Errico
  Presutti.
\newblock Geometry of contours and {P}eierls estimates in {$d=1$} {I}sing
  models with long range interactions.
\newblock {\em J. Math. Phys.}, 46(5):053305, 22, 2005.

\bibitem[Cha03]{Cha03}
D.~Chafa\"\i.
\newblock Glauber versus {K}awasaki for spectral gap and logarithmic {S}obolev
  inequalities of some unbounded conservative spin systems.
\newblock {\em Markov Process. Related Fields}, 9(3):341--362, 2003.

\bibitem[CMR02]{CMR02}
N.~Cancrini, F.~Martinelli, and C.~Roberto.
\newblock The logarithmic {S}obolev constant of {K}awasaki dynamics under a
  mixing condition revisited.
\newblock {\em Ann. Inst. H. Poincar\'e Probab. Statist.}, 38(4):385--436,
  2002.

\bibitem[DF16]{DuoFat16}
M.~H. Duong and M.~Fathi.
\newblock The two-scale approach to hydrodynamic limits for non-reversible
  dynamics.
\newblock {\em Markov Process. Related Fields}, 22(1):1--36, 2016.

\bibitem[Dys69]{Dys69}
Freeman~J. Dyson.
\newblock Existence of a phase-transition in a one-dimensional {I}sing
  ferromagnet.
\newblock {\em Comm. Math. Phys.}, 12(2):91--107, 1969.

\bibitem[FS82]{FrSp82}
J\"urg Fr\"ohlich and Thomas Spencer.
\newblock The phase transition in the one-dimensional {I}sing model with
  {$1/r^{2}$}\ interaction energy.
\newblock {\em Comm. Math. Phys.}, 84(1):87--101, 1982.

\bibitem[GOVW09]{GrOtViWe09}
Natalie Grunewald, Felix Otto, C\'edric Villani, and Maria~G. Westdickenberg.
\newblock A two-scale approach to logarithmic {S}obolev inequalities and the
  hydrodynamic limit.
\newblock {\em Ann. Inst. H. Poincar\'e Probab. Statist.}, 45(2):302--351,
  2009.

\bibitem[GPV88]{GuPaVa88}
M.~Z. Guo, G.~C. Papanicolaou, and S.~R.~S. Varadhan.
\newblock Nonlinear diffusion limit for a system with nearest neighbor
  interactions.
\newblock {\em Comm. Math. Phys.}, 118(1):31--59, 1988.

\bibitem[Gro75]{Gro75}
L.~Gross.
\newblock Logarithmic {S}obolev inequalities.
\newblock {\em Amer. J. Math.}, 97:1061--1083, 1975.

\bibitem[HM16]{HeMe16}
Christopher Henderson and Georg Menz.
\newblock Equivalence of a mixing condition and the {LSI} in spin systems with
  infinite range interaction.
\newblock {\em Stochastic Process. Appl.}, 126(10):2877--2912, 2016.

\bibitem[HS87]{HolStr87}
Richard Holley and Daniel Stroock.
\newblock Logarithmic {S}obolev inequalities and stochastic {I}sing models.
\newblock {\em Journal of Statistical Physics}, 46(5-6):1159--1194, 1987.

\bibitem[Imb82]{Imb82}
John~Z. Imbrie.
\newblock Decay of correlations in the one-dimensional {I}sing model with
  {$J_{ij}=\mid i-j\mid ^{-2}$}.
\newblock {\em Comm. Math. Phys.}, 85(4):491--515, 1982.

\bibitem[KL99]{KipLan99}
Claude Kipnis and Claudio Landim.
\newblock {\em Scaling limits of interacting particle systems}, volume 320 of
  {\em Grundlehren der Mathematischen Wissenschaften [Fundamental Principles of
  Mathematical Sciences]}.
\newblock Springer-Verlag, Berlin, 1999.

\bibitem[KM18]{KwMe17}
Younghak Kwon and Georg Menz.
\newblock Strict convexity of the free energy of the canonical ensemble under
  decay of correlations.
\newblock {\em Journal of Statistical Physics}, Jun 2018.

\bibitem[KM19]{KwMe18a}
Younghak Kwon and Georg Menz.
\newblock Decay of correlations and uniqueness of the infinite-volume gibbs
  measure of the canonical ensemble of 1d-lattice systems.
\newblock {\em Journal of Statistical Physics}, 176(4):836--872, Aug 2019.

\bibitem[Led01a]{Led01}
M.~Ledoux.
\newblock Logarithmic {S}obolev inequalities for unbounded spin systems
  revisted.
\newblock {\em Sem. Probab. XXXV, Lecture Notes in Math, Springer- Verlag},
  1755:167--194, 2001.

\bibitem[Led01b]{Led01a}
Michel Ledoux.
\newblock {\em The concentration of measure phenomenon}, volume~89 of {\em
  Mathematical Surveys and Monographs}.
\newblock American Mathematical Society, Providence, RI, 2001.

\bibitem[LPY02]{LaPaYa02}
C.~Landim, G.~Panizo, and H.~T. Yau.
\newblock Spectral gap and logarithmic {S}obolev inequality for unbounded
  conservative spin systems.
\newblock {\em Ann. Inst. H. Poincar\'e Probab. Statist.}, 38(5):739--777,
  2002.

\bibitem[LY93]{LuYa93}
Sheng~Lin Lu and H.~T. Yau.
\newblock Spectral gap and logarithmic {S}obolev inequality for {K}awasaki and
  {G}lauber dynamics.
\newblock {\em Comm. Math. Phys.}, 156(2):399--433, 1993.

\bibitem[Men11]{Me11}
Georg Menz.
\newblock L{SI} for {K}awasaki dynamics with weak interaction.
\newblock {\em Comm. Math. Phys.}, 307(3):817--860, 2011.

\bibitem[Men14]{Me13}
Georg Menz.
\newblock The approach of {O}tto-{R}eznikoff revisited.
\newblock {\em Electron. J. Probab.}, 19:no. 107, 27, 2014.

\bibitem[MN14]{MeNi14}
Georg Menz and Robin Nittka.
\newblock Decay of correlations in 1{D} lattice systems of continuous spins and
  long-range interaction.
\newblock {\em Journal of Statistical Physics}, 156(2):239--267, 2014.

\bibitem[MO13]{MeOt13}
Georg Menz and Felix Otto.
\newblock Uniform logarithmic {S}obolev inequalities for conservative spin
  systems with super-quadratic single-site potential.
\newblock {\em Ann. Probab.}, 41(3B):2182--2224, 2013.

\bibitem[OR07]{OttRez07}
Felix Otto and Maria~G. Reznikoff.
\newblock A new criterion for the logarithmic {S}obolev inequality and two
  applications.
\newblock {\em Journal of Functional Analysis}, 243(1):121 -- 157, 2007.

\bibitem[Roy99]{Roy99}
G.~Royer.
\newblock Une initiation aux in\'egalit\'es de {S}obolev logarithmiques.
\newblock {\em Cours Sp\'ec., Soc. Math. France}, 1999.

\bibitem[Yau91]{Yau91}
H.~T. Yau.
\newblock Relative entropy and hydrodynamics of {G}inzburg-{L}andau models.
\newblock {\em Lett. Math. Phys.}, 22(1):63--80, 1991.

\bibitem[Yau96]{Yau96}
H.~T. Yau.
\newblock Logarithmic {S}obolev inequality for lattice gases with mixing
  conditions.
\newblock {\em Comm. Math. Phys.}, 181(2):367--408, 1996.

\bibitem[Yos01]{Yos01}
Nobuo Yoshida.
\newblock The equivalence of the log-{S}obolev inequality and a mixing
  condition for unbounded spin systems on the lattice.
\newblock {\em Ann. Inst. H. Poincar\'e Probab. Statist.}, 37(2):223--243,
  2001.

\bibitem[Yos03]{Yos03}
N.~Yoshida.
\newblock Phase transition from the viewpoint of relaxation phenomena.
\newblock {\em Reviews in Mathematical Physics}, 15(07):765--788, 2003.

\bibitem[Zeg96]{Zeg96}
Boguslaw Zegarlinski.
\newblock The strong decay to equilibrium for the stochastic dynamics of
  unbounded spin systems on a lattice.
\newblock {\em Comm. Math. Phys.}, 175(2):401--432, 1996.

\end{thebibliography}

\end{document}